\documentclass{article}
\usepackage{amsmath, amssymb, latexsym}
\usepackage{amsopn,amsfonts, amsthm}
\usepackage[T1]{fontenc}

\newtheorem{theorem}{Theorem}[section]

\newtheorem{corollary}[theorem]{Corollary}

\newtheorem{proposition}[theorem]{Proposition}
\newtheorem{lemma}[theorem]{Lemma}
\theoremstyle{definition}

\newtheorem{definition}[theorem]{Definition}

\newtheorem{question}[theorem]{Question}
\newtheorem{remark}[theorem]{Remark}

\begin{document}

\title{On densely complete metric spaces and extensions of uniformly
continuous functions in $\mathbf{ZF}$}
\author{Kyriakos Keremedis, Eliza Wajch}
\maketitle

\begin{abstract}
A metric space $\mathbf{X}$ is called densely complete if there exists a
dense set $D$ in $\mathbf{X}$ such that every Cauchy sequence of points of $D
$ converges in $\mathbf{X}$. One of the main aims of this work is to prove
that the countable axiom of choice, $\mathbf{CAC}$ for abbreviation, is
equivalent with the following statements:\smallskip

(i) Every densely complete (connected) metric space $\mathbf{X}$ is
complete.\smallskip\ 

(ii) For every pair of metric spaces $\mathbf{X}$ and $\mathbf{Y}$, if $%
\mathbf{Y}$ is complete and $\mathbf{S}$ is a dense subspace of $\mathbf{X}$%
, while $f:\mathbf{S}\rightarrow \mathbf{Y}$ is a uniformly continuous
function, then there exists a uniformly continuous extension $F:\mathbf{X}\to%
\mathbf{Y}$ of $f$.\smallskip

(iii) Complete subspaces of metric spaces have complete closures.\smallskip

(iv) Complete subspaces of metric spaces are closed.\smallskip

It is also shown that the restriction of (i) to subsets of the real line is
equivalent to the restriction $\mathbf{CAC}(\mathbb{R})$ of $\mathbf{CAC}$
to subsets of $\mathbb{R}$. However, the restriction of (ii) to subsets of $%
\mathbb{R}$ is strictly weaker than $\mathbf{CAC}(\mathbb{R})$ because it is
equivalent with the statement that $\mathbb{R}$ is sequential. Moreover,
among other relevant results, it is proved that, for every positive integer $%
n$, the space $\mathbb{R}^n$ is sequential if and only if $\mathbb{R}$ is
sequential. It is also shown that $\mathbb{R}\times\mathbb{Q}$ is not
densely complete if and only if $\mathbf{CAC}(\mathbb{R})$ holds. \bigskip\ 
\newline

\end{abstract}

\section{Notation and terminology}

Let $\mathbf{X}=\langle X,d\rangle$ be a metric space. For $x\in X$ and $%
\varepsilon >0$, let $B_d(x,\varepsilon )=\{y\in X: d(x,y)<\varepsilon \}$.
The topology on $X$ induced by $d$ is denoted by $\tau(d)$. If $A\subseteq X$%
, the diameter of $A$ in $\mathbf{X}$ is denoted by $\delta(A)$ or by $%
\delta_d(A)$ if it is better to point at the metric in use. In the sequel,
boldface letters will denote metric spaces and lightface letters will denote
their underlying sets. If $Y\subseteq X$, then $\mathbf{Y}$ denotes the
metric subspace of $\mathbf{X}$, that is $\mathbf{Y}=(Y, d_{\vert Y})$ where 
$d_{\vert Y}(x,y)=d(x,y)$ for all $x,y\in Y$.\smallskip

If $\tau$ is a topology on a set $Z$ and $Y$ is a subset of $Z$, then, for $%
\tau_{\vert Y}=\{U\cap Y: U\in\tau\}$, the topological space $\langle Y,
\tau_{\vert Y}\rangle$ can be denoted by $Y$ and called the \textit{%
(topological) subspace} of $\langle Z, \tau\rangle$. Topological subspaces
of the metric space $\mathbf{X}$ are topological subspaces of the
topological space $\langle X, \tau(d)\rangle$. \smallskip

The set of all finite von Neumann ordinal numbers is denoted by $\omega$,
while $\mathbb{N}=\omega\setminus\{0\}$. We recall that if $n\in\omega$,
then $n+1=n\cup\{n\}$.\smallskip The power set of $X$ is denoted by $%
\mathcal{P}(X)$.\smallskip

In Definitions \ref{d1.0}-\ref{d1.3}, we recall several known notions.

\begin{definition}
\label{d1.0} Let $A\subseteq X$. We denote by $A^{\ast}$ the set of all
points $x\in X$ which have the property that there exists a sequence $%
(x_{n})_{n\in \mathbb{N}}$ of points of $A\setminus\{x\}$ such that $%
\lim\limits_{n\to\infty}d(x_{n}, x)=0$. Then the set $\widetilde{A}=A\cup
A^{\ast}$ is called the \textit{sequential closure} of $A$ in $\mathbf{X}$.
The set $A$ is called \textit{sequentially closed} in $\mathbf{X}$ if $%
A^{\ast}\subseteq A$.
\end{definition}

\begin{definition}
\label{d1.2} The space $\mathbf{X}$ is called:

\begin{enumerate}
\item[(i)] \textit{sequential} if each sequentially closed subset of $%
\mathbf{X}$ is closed in $\mathbf{X}$;

\item[(ii)] \textit{Fr\'{e}chet-Urysohn} if for every $F\in \mathcal{P}(X)$,
all accumulation points of $F$ belong to the sequential closure of $F$;

\item[(iii)] \textit{complete} if each of its Cauchy sequences converges in $%
\mathbf{X}$;

\item[(iv)] \textit{Cantor complete} if, for every descending family $%
\{G_{n}: n\in \mathbb{N}\}$ of non-empty closed subsets of $\mathbf{X}$ with 
$\lim\limits_{n\rightarrow \infty }\delta(G_{n})=0$, the set $%
\bigcap\limits_{n\in\mathbb{N}}G_{n}$ is a singleton.
\end{enumerate}
\end{definition}

\begin{definition}
\label{d1.3}  Let $f$ be a function from $\mathbf{X}$ into a metric space $%
\mathbf{Y}$. Then $f$ is called:

\begin{enumerate}
\item[(i)] \textit{sequentially continuous at} $x_{0}\in X$ if for every
sequence $(x_{n})_{n\in \mathbb{N}}$ of points of $X$ such that $%
(x_{n})_{n\in \mathbb{N}}$ converges in $\mathbf{X}$ to $x_0$, the sequence $%
(f(x_n))_{n\in\mathbb{N}}$ converges in $\mathbf{Y}$ to $f(x_0)$;

\item[(ii)] \textit{sequentially continuous} if $f$ is sequentially
continuous at each point of $\mathbf{X}$.
\end{enumerate}
\end{definition}

To a great extent, this article is about new concepts introduced in the
following definition:

\begin{definition}
\label{d1.1} (i ) If $\rho$ is a metric on a set $Y$ and $\mathbf{Y}=\langle
Y, \rho\rangle$, then:

\begin{enumerate}
\item[(a)] for $D\subseteq Y$, the metric space $\mathbf{Y}=\langle Y,
\rho\rangle$ is called \textit{densely complete with respect to} $D$ if $D$
is dense in $\mathbf{Y}$ and every Cauchy sequence of points of $D$
converges in $\mathbf{Y}$;

\item[(b)] the metric space $\mathbf{Y}$ is called \textit{densely complete}
if there exists a subset $D$ of $Y$ such that $\mathbf{Y}$ is densely
complete with respect to $D$;

\item[(c)] the metric $\rho$ is called \textit{densely complete} if the
metric space $\mathbf{Y}$ is densely complete.
\end{enumerate}

(ii) A topological space $\langle Y, \tau\rangle$ is called \textit{densely
completely metrizable} if there exists a densely complete metric $\rho$ on $Y
$ such that $\tau=\tau(\rho)$.
\end{definition}

The symbol $\mathbb{R}$ denotes the set of real numbers equipped with the
metric induced by the standard absolute value. As usual, for $n\in\mathbb{N}$
and $x\in\mathbb{R}^n$, let $\Vert x\Vert=\sqrt{\sum_{i\in n}x(i)^2}$
and, for $x, y\in\mathbb{R}^n$, let $d_{e(n)}(x, y)=\Vert x-y\Vert$%
. Then $d_{e(n)}$ is the Euclidean metric on $\mathbb{R}^n$ and, if it is
not stated otherwise, $\mathbb{R}^n$ will denote both the metric space $%
\langle \mathbb{R}^n, d_{e(n)}\rangle$ and the set $\mathbb{R}^n$, as well
as the topological space $\langle\mathbb{R}^n, \tau(d_{e(n})\rangle$%
.\smallskip

We recall that if $\{\mathbf{X}_{n}: n\in \mathbb{N}\}$ is a family of
metric spaces $\mathbf{X}_{n}=\langle X_{n},d_{n}\rangle$, while $%
X=\prod\limits_{n\in \mathbb{N}}X_{n}$ and, for all $n\in \mathbb{N}$ and $%
a,b\in X_{n}$, $\rho _{n}(a,b)=\min \{1,d_{n}(a,b)\}$, then the function $%
d:X\times X\rightarrow \mathbb{R}$ given by: 
\begin{equation}
d(x,y)=\sum\limits_{n\in \mathbb{N}}\frac{\rho (x_{n},y_{n})}{2^{n}}
\label{2}
\end{equation}%
is a metric such that the topology $\tau (d)$ coincides with the product
topology of the family of topological spaces $\{\langle X_{n},\tau
(d_{n})\rangle: n\in \mathbb{N}\}$. In the sequel, we shall always assume
that whenever a family $\{\mathbf{X}_{n}: n\in \mathbb{N}\}$ of metric
spaces $\mathbf{X}_{n}=\langle X_{n},d_{n}\rangle$ is given, then $%
X=\prod\limits_{n\in \mathbb{N}}X_{n}$ and $\mathbf{X}=\langle X,d\rangle$
where $d$ is the metric on $X$ defined by (\ref{2}).\smallskip

A set which contains a countably infinite subset is called \textit{%
Dedekind-infinite}. That $X$ is Dedekind-infinite is denoted by $\mathbf{DI}%
(X)$. A set which is not Dedekind-infinite is called \textit{Dedekind-finite}%
. By universal quantifying over $X$, $\mathbf{DI}(X)$ gives rise to the
choice principle $\mathbf{IDI}$: $\forall X(X$ infinite $\rightarrow $ $%
\mathbf{DI}(X))$ that is, \textquotedblleft every infinite set is
Dedekind-infinite\textquotedblright\ (Form 9 of \cite{hr}, Definition
2.13(1) in \cite{herl}). Below we list some other weak forms of the axiom of
choice we shall deal with in the sequel. For the known forms, we quote in
their statements, the form number under which they are recorded in \cite{hr}
or we refer to their definitions in \cite{herl}. The rest of the forms are
new here. We recall that if $\mathcal{A}$ is an infinite collection of sets,
every element of the product $\prod\limits_{A\in\mathcal{A}}A$ is called a 
\textit{choice function} of $\mathcal{A}$, while every choice function of an
infinite subcollection of $\mathcal{A}$ is called a \textit{partial choice
function} of $\mathcal{A}$. Given a pairwise disjoint family $\mathcal{A}$,
a set $C$ is called a partial choice set of $\mathcal{A}$ if the set $\{A\in 
\mathcal{A}: C\cap A \text{ is a singleton }\}$ is infinite.

\begin{itemize}
\item $\mathbf{CAC}(X)$: Every non-empty countable family of non-empty
subsets of $X$ has a choice function.

\item $\mathbf{CAC}$ (Form 8 in \cite{hr}, Definition 2.5 in \cite{herl}):
For every infinite set $X$, $\mathbf{CAC}(X)$. Equivalently, every countable
family of pairwise disjoint non-empty sets has a partial choice set (cf.
Form [8 B] in \cite{hr}).

\item $\mathbf{CMC}$ (Form 126 in \cite{hr}, Definition 2.10 in \cite{herl}%
): For every collection $\{A_n: n\in\omega\}$ of non-empty sets there exists
a collection $\{F_n: n\in\omega\}$ of non-empty finite sets such that $%
F_n\subseteq A_n$ for each $n\in\omega$.

\item $\mathbf{CAC}_{fin}$ (Form 10 in \cite{hr}, Definition 2.9(3) in \cite%
{herl}): $\mathbf{CAC}$ restricted to countable families of non-empty finite
sets.

\item $\mathbf{CAC(\mathbb{R})}$ (Form 94 in \cite{hr}, Definition 2.9(1) in 
\cite{herl}): $\mathbf{CAC}$ restricted to families of non-empty subsets of $%
\mathbb{R}$. Equivalently, every countable family of pairwise disjoint
non-empty subsets of the real line has a partial choice set (see, e.g., \cite%
{hkrt}).

\item $\omega -\mathbf{CAC}(\mathbb{R})$ (Definition 4.56 in \cite{herl}):
For every family $\mathcal{A}=\{A_{i}: i\in \omega \}$ of non-empty subsets
of $\mathbb{R}$ there exists a family $\mathcal{B}=\{B_{i}: i\in \omega \}$
of countable non-empty subsets of $\mathbb{R}$ such that for each $i\in
\omega $, $B_{i}\subseteq A_{i}$.

\item $\mathbf{IDI(\mathbb{R})}$ (Form 13 in \cite{hr}, Definition 2.13(2)
in \cite{herl}): $\mathbf{IDI}$ restricted to subsets of $\mathbb{R}$.

\item $\mathbf{DC}$ (Form 43 in \cite{hr}, Definition 2.11 (1) in \cite{herl}%
): For every pair $\langle X, \rho\rangle$, where $X$ and $\rho\subseteq
X\times X$ is such that for each $x\in X$ there exists $y\in X$ with $%
\langle x,y\rangle\in\rho$, then there exists a sequence $(x_n)_{n\in\mathbb{%
N}}$ of points of $X$ such that $\langle x_n, x_{n+1}\rangle\in\rho$ for
each $n\in\mathbb{N}$.

\item $\mathbf{BPI}$ (Form 14 in \cite{hr}, Definition 2.15 in \cite{herl}):
Every Boolean algebra has a prime ideal.
\end{itemize}

It is known that $\mathbf{BPI}$ is equivalent to the statement: The
Tychonoff product of compact $T_{2}$\ spaces is compact (see, e.g., \cite{rs}%
, Form [14 J] in \cite{hr} and Theorem 4.70 in \cite{herl}).

Our main aim is to investigate the following newly defined forms:

\begin{itemize}
\item $\mathbf{DCC}(\mathbf{X})$: Every densely complete subspace of $%
\mathbf{X}$ is complete.

\item $\mathbf{DCC}$: Every densely complete metric space is complete.

\item $\mathbf{DCC}_{c}(\mathbf{X})$: Every densely complete connected
subspace of $\mathbf{X}$ is complete.

\item $\mathbf{DCC}_c$: Every densely complete connected metric space is
complete.

\item $\mathbf{UCE}(\mathbf{X})$: For every metric subspace $\mathbf{Z}$ of $%
\mathbf{X}$ and every complete metric subspace $\mathbf{Y}$ of $\mathbf{X}$,
if $S$ is a dense subset of $\mathbf{Z}$ and $f:\mathbf{S}\rightarrow 
\mathbf{Y}$ is uniformly continuous, then there exists a uniformly
continuous extension $F:\mathbf{Z}\to\mathbf{Y}$ of $f$.\smallskip

\item $\mathbf{UCE}$: For every metric space $\mathbf{X}$, $\mathbf{UCE}(%
\mathbf{X})$ holds.\smallskip

\item $\mathbf{UCE}(\mathbf{X},\mathbb{R})$: For every metric subspace $%
\mathbf{Z}$ of $\mathbf{X} $ and every complete metric subspace $\mathbf{Y}$
of $\mathbb{R}$, if $S$ is a dense subset of $\mathbf{Z}$ and $f:\mathbf{S}%
\rightarrow \mathbf{Y}$ is uniformly continuous, then there exists a
uniformly continuous extension $F:\mathbf{Z}\to\mathbf{Y}$ of $f$.\smallskip

\item $\mathbf{UCE}_{c}$: For every compact metric space $\mathbf{X}$ and
every complete metric space $\mathbf{Y}$, if $S$ is a dense subset of $%
\mathbf{X}$ and $f:\mathbf{S}\rightarrow \mathbf{Y}$ is a uniformly
continuous function, then there exists a uniformly continuous extension $F:%
\mathbf{X}\to\mathbf{Y}$ of $f$.

\item $\mathbf{UCE}_{cc}(\mathbf{X})$: For every compact and connected
metric subspace $\mathbf{Z} $ of $\mathbf{X}$ and every complete metric
subspace $\mathbf{Y}$ of $\mathbf{X}$, if $S$ is a dense subset of $\mathbf{Z%
}$ and $f:\mathbf{S}\rightarrow \mathbf{Y}$ is uniformly continuous, then
there exists a uniformly continuous extension $F:\mathbf{Z}\rightarrow 
\mathbf{Y}$ of $f$.

\item $\mathbf{UCE}_{cc}$: For every metric space $\mathbf{X}$, $\mathbf{UCE}%
_{cc}(\mathbf{X})$ holds.

\item $\mathbf{UCE}_{cc}(\mathbf{X},\mathbb{R})$: For every compact and
connected metric subspace $\mathbf{Z}$ of $\mathbf{X}$ and every complete
metric subspace $\mathbf{Y}$ of $\mathbb{R}$, if $S$ is a dense subset of $%
\mathbf{Z}$ and $f:\mathbf{S}\rightarrow \mathbf{Y}$ is uniformly
continuous, then there exists a uniformly continuous extension $F:\mathbf{Z}%
\rightarrow \mathbf{Y}$ of $f$.\smallskip

\item $\mathbf{UCE}_{c}(\mathbb{R})$: For every connected metric subspace $%
\mathbf{X}$ of $\mathbb{R}$ and every complete metric subspace $\mathbf{Y}$
of $\mathbb{R} $, if $S$ is a dense subset of $\mathbf{X}$ and $f:\mathbf{S}%
\rightarrow \mathbf{Y}$ is a uniformly continuous function, then there
exists a uniformly continuous extension $F:\mathbf{X}\to\mathbf{Y}$ of $f$.
\end{itemize}

Usually, a topological space which is simultaneously compact and connected
is called a continuum. Since every non-empty continuum of $\mathbb{R}$ is
either a singleton or an interval $[a,b]$ where $a,b\in \mathbb{R}$ and $a<b$%
, then the following proposition holds:

\begin{proposition}
\label{p1.1} $\mathbf{UCE}_{cc}(\mathbb{R}, \mathbb{R})$ and $\mathbf{UCE}%
_{cc}(\mathbb{R})$ are equivalent.
\end{proposition}

We recall that if $\mathbf{X}, \mathbf{Y}$ are metric spaces, $S$ is a dense
set of $\mathbf{X}$, while $f:\mathbf{S}\to\mathbf{Y}$ is extendable to a
uniformly continuous mapping from $\mathbf{X}$ to $\mathbf{Y}$, then there
exists a unique uniformly continuous extension of $f$ over $X$.

\section{Introduction and some preliminary results}

In this paper, the intended context for reasoning and statements of theorems
will be the Zermelo-Fraenkel set theory $\mathbf{ZF}$ without the axiom of
choice $\mathbf{AC}$. As usual, the system $\mathbf{ZF+AC}$ is denoted by $%
\mathbf{ZFC}$. Our intention here is to study the set-theoretic strength of
the new statements $\mathbf{DCC}$, $\mathbf{DCC}_{c}$, $\mathbf{DCC}(\mathbb{%
R})$, $\mathbf{DCC}_{c}(\mathbb{R})$, $\mathbf{UCE}$ and $\mathbf{UCE}(%
\mathbb{R})$, as well as to sort them in the hierarchy of choice principles.
We also solve some non-trivial open problems posed in \cite{OPWZ}.

The paper is organized as follows. In Section 2, we give some preliminary
results. Section 3 concerns $\mathbf{UCE}(\mathbb{R})$, $\mathbf{UCE}_c(%
\mathbb{R})$, $\mathbf{DCC}(\mathbb{R})$, $\mathbf{DCC}_c(\mathbb{R})$, as
well as some modifications of $\mathbf{UCE}(\mathbb{R})$. In Section 3,
among other results, it is proved that $\mathbf{UCE}(\mathbb{R})$ and $%
\mathbf{UCE}_c(\mathbb{R})$ are both equivalent to the sentence: $\mathbb{R}$
is sequential. It is also proved in Section 3 that $\mathbf{DCC}(\mathbb{R})$%
, $\mathbf{DCC}_c(\mathbb{R})$ and $\mathbf{CAC}(\mathbb{R})$ are all
equivalent. Moreover, it is proved in Section 3 that, for each $n\in\mathbb{N%
}$, the space ${\mathbb{R}}^n$ is sequential if and only if $\mathbb{R}$ is
sequential. This leads to immediate solutions of Problems 6.6 and 6.11 posed
in \cite{OPWZ}.

In Section 4, it is proved that $\mathbf{DCC}$, $\mathbf{DCC}_{c}$ and $%
\mathbf{UCE}$ are all equivalent with $\mathbf{CAC}$; however, $\mathbf{UCE}%
_{c}$ does not imply $\mathbf{CAC}$ in $\mathbf{ZF}^{0}$. It has been
established in \cite{kker} that $\mathbf{UCE}$ implies $\mathbf{CAC}_{fin}$.
Therefore, in the forthcoming Theorem \ref{t4.12}, we improve the latter
result by showing that $\mathbf{CAC}$ and $\mathbf{UCE}$ are actually
equivalent. Moreover, we investigate the problem on when countable products
of densely complete metric spaces can be densely complete. We also show that 
$\mathbf{CAC}(\mathbb{R})$ is equivalent with the statement that $\mathbb{R}%
\times\mathbb{Q}$ is not densely complete. \smallskip

It is well-known that, in the absence of the axiom of choice $\mathbf{AC}$,
fundamental results of elementary analysis and topology may fail severely in
some $\mathbf{ZF}$ models. As an example of this kind of models, we mention
here the famous Cohen Basic Model $\mathcal{M}$1 in \cite{hr}. It is known
that in $\mathcal{M}$1, the set $A$ of all added Cohen reals is an infinite,
Dedekind-finite subset of $\mathbb{R}$ disjoint from $\mathbb{Q}$. It is
known that the following disasters may strike in $\mathcal{M}$1:\smallskip

(i) $A$ is a complete subset of $\mathbb{R}$ with $\overline{A}=\mathbb{R}%
\neq A$. Hence, complete metric subspaces of $\mathbb{R}$ need not be closed
in $\mathcal{M}$1.

(ii) $A=\widetilde{A}\neq \overline{A}=\mathbb{R}$. So, $\mathbb{R}$ is not
sequential in $\mathcal{M}$1.

(iii) $0\in \overline{A}\backslash A$. This means that $\mathbb{R}$ is not Fr%
\'{e}chet-Urysohn in $\mathcal{M}$1.

(iv) Any function $f:\mathbf{A}\rightarrow \mathbb{R}$ is sequentially
continuous. Thus, sequential continuity of real-valued functions defined on
complete metric subspaces of $\mathbb{R}$ is not equivalent to the usual $%
\varepsilon -\delta $ definition of continuity in $\mathcal{M}$1.

(v) $\mathbf{A}$ is not separable. Therefore, $\mathbb{R}$ is not
hereditarily separable in $\mathcal{M}$1.

(vi) The open cover $\{(a-1,a+1): a\in A\}$ of $\mathbb{R}$ has no countable
subcover. So, $\mathbb{R}$ is not Lindel\"{o}f in $\mathcal{M}$1.

(vii) $(0,1)\cap A$ is a complete, non-closed metric subspace of the
connected subspace $(0,1) $ of $\mathbb{R}$, but its closure in $(0,1)$ is
not complete. So, closures of complete metric subspaces of connected metric
spaces need not be complete in $\mathcal{M}$1!\smallskip

For more disasters of this kind, we refer the reader to \cite{herl}.

\begin{proposition}
\label{p2.1} $\mathbf{DCC}_c(\mathbb{R})$ implies $\mathbf{IDI(\mathbb{R})}$.
\end{proposition}

\begin{proof}
Assume the contrary and let $D$ be an infinite Dedekind-finite subset of $%
\mathbb{R}$. N. Brunner showed in \cite{br} that if there is a Dedekind-finite set $D$ of $\mathbb{R}$, then there exists a dense one also. This is why we
may assume that $D$ is dense in $\mathbb{R}$. Let $H=D\cap (0, 1)$. Since $H$
contains no countably infinite subsets, it follows that every Cauchy
sequence of points of $H$ is eventually constant, hence convergent in $H$,
so in $(0, 1)$. Hence $(0, 1)$ is a densely complete connected, not complete metric
subspace of $\mathbb{R}$.
\end{proof}

The proof to Proposition \ref{p2.1} shows that each one of the statements $%
\mathbf{DCC}_c(\mathbb{R})$, $\mathbf{DCC}(\mathbb{R})$, $\mathbf{DCC}_c$
and $\mathbf{DCC}$ fails in the model $\mathcal{M}$1. In consequence, we
have the following corollary:

\begin{corollary}
\label{p2.2} None of $\mathbf{DCC}_c(\mathbb{R})$, $\mathbf{DCC}(\mathbb{R})$%
, $\mathbf{DCC}_c$, $\mathbf{DCC}$ is a theorem of $\mathbf{ZF}$.
\end{corollary}

Let us list in the following theorem some known results we will need in the
sequel.

\begin{theorem}
\label{t2.3} (i) (Cf. \cite{gg2}, \cite{nak} and pp. 74--75 of \cite{herl}.)
The following are equivalent:\newline
(a) $\mathbf{CAC}(\mathbf{\mathbb{R}})$;\newline
(b) every subspace of $\mathbb{R}$ is sequential;\newline
(c) $\mathbb{R}$ is Fr\'{e}chet-Urysohn;\newline
(d) for every $x\in \mathbb{R}$ and every function $f:\mathbb{R}\rightarrow 
\mathbb{R}$, if $f$ is sequentially continuous at $x$, then $f$ is
continuous at $x$.\newline
(ii) (Cf. \cite{gg2}.) The following are equivalent:\newline
(a) $\mathbf{CAC}$;\newline
(b) every metric space is sequential;\newline
(c) every metric space is\smallskip\ Fr\'{e}chet-Urysohn;\newline
(d) every complete metric space is Cantor complete.\newline
(iii)(a) (Cf. \cite{herl}, p.75.) $\mathbb{R}$ is sequential if and only if
every complete metric subspace of $\mathbb{R}$ is closed.\newline
(b) (Cf. \cite{herl}, p.76.) $\omega -\mathbf{CAC}(\mathbf{\mathbb{R}})$
implies $\mathbb{R}$ is sequential.\newline
(c) (Cf. \cite{gg1}.) If there is a complete, non-closed metric subspace of $%
\mathbb{R}$, then there exists a dense, complete, non-closed metric subspace
of $\mathbb{R}$.\newline
(iv) (a) (Cf. \cite{l}.) The family of all non-empty closed subsets of the
real line has a choice function.\newline
(b) (Cf. \cite{kt}.) If the family of all non-empty closed subsets of a
compact metric space has a choice set, then the space is separable.
\end{theorem}

Clearly, every complete metric space is trivially densely complete, and the
proposition: \textit{every densely complete metric space is complete} is a
well-known criterion of completeness of metric spaces in $\mathbf{ZF}+%
\mathbf{CAC}$. It is also known that the Baire Category Theorem can fail for
a complete metric space in a model of $\mathbf{ZF}$ and that every separable
completely metrizable space is a Baire space (see, e.g., Section 4.10 of 
\cite{herl}). For the reader's convenience we supply a proof to the
following more general theorem here:

\begin{theorem}
\label{t2.4} Let $\mathbf{X}=\langle X,d\rangle$ be a densely complete with
respect to a set $D$ metric space. Then the following conditions are
satisfied:\newline
(i) If $\mathbf{CAC}(D)$ holds, then $\mathbf{X}$ is complete.\newline
(ii) If $D$ is well-orderable, then $\mathbf{X}$ is a Baire space.
\end{theorem}

\begin{proof}
(i) Assume that $\mathbf{CAC}(D)$ holds. To show that $\mathbf{X}$ is
complete, fix a Cauchy sequence $(x_{n})_{n\in \mathbb{N}}$ of $\mathbf{X}$
and, for each $n\in \mathbb{N}$, put $D_{n}=B_d(x_{n},\frac{1}{n})\cap D$.
Since $D$ is dense in $\mathbf{X}$, the sets $D_{n}$ are non-empty. In view
of $\mathbf{CAC}(D)$, there is a sequence $(y_{n})_{n\in \mathbb{N}}$ such
that $y_{n}\in D_{n}$ for each $n\in \mathbb{N}$. We claim that $%
(y_{n})_{n\in \mathbb{N}}$ is Cauchy. To see this, fix $\varepsilon >0$.
Since $(x_{n})_{n\in \mathbb{N}}$ is Cauchy, it follows that there exists $%
n_{0}\in \mathbb{N}$ such that $1/n_{0}<\varepsilon /3$ and $d(x_{n},x_{m})<\varepsilon /3$ whenever $n,m\geq
n_{0}$. For all natural numbers $n,m\geq n_{0}$, we have: 
\begin{equation*}
d(y_{n},y_{m})\leq d(y_{n},x_{n})+d(x_{n},x_{m})+d(x_{m},y_{m})\leq
1/n+\varepsilon /3+1/m<\varepsilon \text{.}
\end{equation*}%
Hence, $(y_{n})_{n\in \mathbb{N}}$ is Cauchy as claimed. Therefore, $%
(y_{n})_{n\in \mathbb{N}}$ converges to some point $x\in X$. Since for all $%
n\in \mathbb{N}$, $d(x_{n},x)\leq d(x_{n},y_{n})+d(y_{n},x)$, we get that 
\begin{equation*}
\lim\limits_{n\to\infty} d(x_{n},x)\leq \lim\limits_{n\to\infty}d(x_{n},y_{n})+\lim\limits_{n\to\infty} d(y_{n},x)=0\text{.}
\end{equation*}%
Thus, $\lim\limits_{n\to\infty}d( x_{n},x)=0$. This is why $\mathbf{X}$ is complete as required.\smallskip

(ii) Let us notice that if $D$ is well-orderable, then the base $\mathcal{B}%
=\{ B_d(x, \frac{1}{n}): x\in D \text{ and } n\in\mathbb{N}\}$ of $\mathbf{X}
$ is well-orderable, so the usual proof in $\mathbf{ZFC}$ that complete
metric spaces are Baire can be easily adopted to get a proof in $\mathbf{ZF}$
that $\mathbf{X}$ is Baire.
\end{proof}

\begin{corollary}
\label{c2.5} (i) $\mathbf{CAC}$ implies $\mathbf{DCC}$.\newline
(ii) $\mathbf{CAC}(\mathbb{R})$ implies $\mathbf{DCC}(\mathbb{R})$.\newline
(iii) For every metric space $\mathbf{X}=\langle X,d\rangle$, $\mathbf{CAC}%
(X)$ implies $\mathbf{DCC}(\mathbf{X})$.
\end{corollary}

\begin{theorem}
\label{t2.d0} $\mathbf{DC}$ is equivalent with the statement: Every densely
complete metric space is Baire.
\end{theorem}

\begin{proof}
It follows from Theorem 4.106 of \cite{herl} that if every densely complete
metric space is Baire, then $\mathbf{DC}$ holds. Now, assume $\mathbf{DC}$
and assume that $\mathbf{X}=\langle X,d\rangle $ is a densely complete
metric space. Let $D\subseteq X$ be such that $\mathbf{X}$ is densely
complete with respect to a set $D$. To show that $\mathbf{X}$ is Baire, we
slightly modify the part that (1) $\leftrightarrow $ (2) of the proof to
Theorem 4.106 in \cite{herl}. Let $(G_{n})_{n\in \omega }$ be a sequence of
open dense sets in $\mathbf{X}$ such that $G_{n+1}\subseteq G_{n}$ for each $%
n\in \omega $. Let $G\in \tau (d)$ and 
\begin{equation*}
Y=\{\langle n,x,r\rangle \in \omega \times D\times \mathbb{R}:0<r<\frac{1}{%
2^{n}}\text{ and }\overline{B_{d}(x,r)}\subseteq G\cap G_{n}\}.
\end{equation*}%
Define a binary relation $\rho $ on $Y$ as follows: if $y=\langle
n,x,r\rangle \in Y$ and $\overline{y}=\langle \overline{n},\overline{x},%
\overline{r}\rangle \in Y$, then: 
\begin{equation*}
\langle y,\overline{y}\rangle \in \rho \leftrightarrow (n<\overline{n}\text{
and }\overline{B_{d}(\overline{x},\overline{r})}\subseteq B_{d}(x,r)).
\end{equation*}%
It follows from $\mathbf{DC}$ that there exists a sequence $(y_{n})_{n\in
\omega }$ of points of $Y$ such that $\langle y_{n},y_{n+1}\rangle \in \rho $
for each $n\in \omega $. If $y_{n}=\langle m_{n},x_{n},r_{n}\rangle $ for
each $n\in \omega $, then $(x_{n})_{n\in \omega }$ is a Cauchy sequence of
points of $D$. Since $\mathbf{X}$ is densely complete with respect to $D$,
the sequence $(x_{n})_{n\in \omega }$ converges in $\mathbf{X}$ to a point $%
x $. Of course, $x\in G\cap \bigcap\limits_{n\in \omega }G_{n}$, so $\mathbf{X}$ is
Baire.
\end{proof}

\begin{theorem}
\label{t2.6} Let $\mathbf{X}=\langle X,d\rangle$\ be a metric space and $%
Y\subseteq X$. Then the following conditions are satisfied:\newline
(i) If $b\in \overline{Y}\backslash Y$, while $Z=Y\cup \{b\}$, then the
identity function $f:\mathbf{Y}\rightarrow \mathbf{Y}$ does not extend
continuously over $Z$.\newline
(ii) If $\mathbf{X}$ is complete, then $\mathbf{UCE}(\mathbf{X})$ implies $%
\mathbf{X}$ is sequential.\newline
\end{theorem}

\begin{proof}
(i) Assume the contrary and let $F:\mathbf{Z}\rightarrow \mathbf{Y}$ be a
continuous extension of $f$ over $Z$. Then $F(b)=a\in Y$. Let $\varepsilon
=d(a,b)>0$ ($a\neq b$), and fix, by the continuity of $F$ at $b,$ a number $%
\delta \in (0,\varepsilon /2)$ such that: 
\begin{equation*}
\text{for each }x\in Z\text{ with }d(x,b)<\delta ,d(F(x),F(b))<\varepsilon /2%
\text{.}
\end{equation*}%
Since $Y$ is dense in $\mathbf{Z}$, we can fix $x\in Y$ with $d(x,b)<\delta $. Since\ 
\begin{equation*}
d(x,a)=d(F(x),F(b))<\varepsilon /2\text{,}
\end{equation*}%
we get, 
\begin{equation*}
\varepsilon =d(a,b)\leq d(b,x)+d(x,a)<\varepsilon /2+\varepsilon
/2=\varepsilon \text{.}
\end{equation*}%
The contradiction obtained implies that $f$ does not extend continuously over $Z$. This finishes
the proof of (i).\smallskip

(ii) Suppose that $\mathbf{X}=\langle X,d\rangle$ is a
complete metric space having a sequentially closed subset $A$ such that $%
\overline{A}\setminus\widetilde{A}\neq \emptyset $. Then the metric subspace $\mathbf{%
A}$ of $\mathbf{X}$ is complete. It follows from (i) that, for every $a\in 
\overline{A}\setminus A$, the identity function $f:A\rightarrow A$ does not
extend continuously over $A\cup \{a\}$. This contradicts $\mathbf{UCE}(\mathbf{X}%
)$ and finishes the proof of (ii).\medskip
\end{proof}

Another well-known result in elementary analysis which fails in the model $%
\mathcal{M}$1 is the proposition $\mathbf{UCE}(\mathbb{R})$. Indeed, the set 
$A$ of all added Cohen reals is complete and dense in the metric subspace $%
\mathbf{X}$ of $\mathbb{R}$ where $X=A\cup \{0\}$. Hence, by Theorem \ref%
{t2.6}(i), the uniformly continuous function $f:\mathbf{A}\rightarrow 
\mathbf{A}$ does not extend continuously over $X$. Therefore, $\mathbf{UCE}(%
\mathbb{R})$ and, in consequence, $\mathbf{UCE}$ fail in $\mathcal{M}$%
1.\smallskip

\begin{corollary}
\label{c2.7} Neither $\mathbf{UCE}$ nor $\mathbf{UCE}(\mathbb{R})$ is a
theorem of $\mathbf{ZF}$.
\end{corollary}

\section{On $\mathbf{UCE}(\mathbb{R})$ and $\mathbf{DCC}(\mathbb{R}$)}

Our first result in this section shows that $\mathbf{UCE}(\mathbb{R})$ is
equivalent to the strictly weaker\ consequence \textquotedblleft $\mathbb{R}$
is sequential\textquotedblright\ of $\mathbf{CAC}(\mathbb{R})$. Before we
proceed, we need to establish the following lemma:

\begin{lemma}
\label{l3.1} Let $\mathbf{Y}=\langle Y,\rho\rangle$ be a Cantor complete
metric space, let $\mathbf{X}=\langle X,d\rangle$ be a metric space and $S$
a dense set in $\mathbf{X}$. Suppose that $f:\mathbf{S}\rightarrow \mathbf{Y}
$ is uniformly continuous. Then there exists a uniformly continuous
extension $F:\mathbf{X}\rightarrow \mathbf{Y}$ of $f$.
\end{lemma}

\begin{proof}
For every $x\in X,$ let%
\begin{equation*}
\mathcal{F}(x)=\{\overline{f[B_d(x,1/n)\cap S]}: n\in \mathbb{N}\}\text{.}
\end{equation*}%
Since $S$ is dense in $\mathbf{X}$, each set from $\mathcal{F}(x)$ is
non-empty and closed in $\mathbf{Y}$. By the uniform continuity of $f$, for
each $k\in \mathbb{N},$ there exists $n_{k}\in \mathbb{N}$ such that $%
n_{k}>k $ and the inclusion $f[B_d(x,1/{n_{k}})\cap S]\subseteq B_{\rho
}(f(x),1/k)$ holds for each $x\in X$. Thus, for each $x\in X$, 
\begin{equation*}
\lim_{n\rightarrow \infty }{\delta}_{\rho}(\overline{f[B_{d}(x,1/n)\cap S]%
})=0\text{.}
\end{equation*}%
Therefore, $\bigcap \mathcal{F}(x)$ cannot have two distinct points.
Moreover, since $\mathbf{Y}$ is Cantor complete, for each $x\in X$, the set $%
\bigcap \mathcal{F}(x)$ is non-empty. Hence, we can define a mapping $F:%
\mathbf{X}\rightarrow \mathbf{Y}$ by putting $F(x) $ as the unique point of $%
\bigcap \mathcal{F}(x)$ for each $x\in X$. As in standard textbooks of
elementary analysis or topology, it is a routine work to verify that $F$ is the required
extension of $f$ (see, e.g., the proof to Lemma 4.3.16 in \cite{En}).
\end{proof}

\begin{corollary}
\label{c3.2} For a given $m\in\mathbb{N}$, suppose that $\mathbf{Y}$ is a
complete, closed metric subspace of $\mathbb{R}^m$. Let $\mathbf{X}=\langle
X,d\rangle$ be a metric space, $S$ a dense set in $\mathbf{X}$ and $f:%
\mathbf{S}\rightarrow \mathbf{Y}$ a uniformly continuous function. Then
there exists a uniformly continuous extension $F:\mathbf{X}\to\mathbf{Y}$ of 
$f$.
\end{corollary}

\begin{proof}
It suffices to apply Lemma \ref{l3.1} and a simple observation that, since
every bounded closed subset of $\mathbb{R}^m$ is compact, every complete
closed metric subspace of $\mathbb{R}^m$ is Cantor complete.
\end{proof}

\begin{theorem}
\label{t3.3}(i) The following are equivalent:

\begin{enumerate}
\item[(a)] $\mathbf{UCE}(\mathbb{R})$;

\item[(b)] $\mathbf{UCE}_{c}(\mathbb{R})$;

\item[(c)] $\mathbf{UCE}_{cc}(\mathbb{R})$;

\item[(d)] $\mathbb{R}$ is sequential.
\end{enumerate}

(ii) It is relatively consistent with $\mathbf{ZF}$, $\mathbf{UCE}(\mathbb{R}%
)$ and the negation of $\mathbf{CAC}(\mathbb{R})$. Hence, $\mathbf{UCE}(%
\mathbb{R})$ is strictly weaker than $\mathbf{CAC}(\mathbb{R})$.\newline
\end{theorem}

\begin{proof}
(i) The implications (a) $\rightarrow $ (b) $\rightarrow $ (c) are obvious.

(c) $\rightarrow $ (d) We assume (c). In view of Theorem \ref{t2.3}(iii)(a),
to prove that $\mathbb{R}$ is sequential, it suffices to show that every
complete metric subspace of $\mathbb{R}$ is closed. Assume the contrary and fix a
complete, non-closed metric subspace $\mathbf{S}$ of $\mathbb{R}$. By Theorem \ref%
{t2.3} (iii) (c), we may assume that $S$ is dense in $\mathbb{R}$. Fix $a\in 
\overline{S}\setminus S$. Then, for $Y=S\cap [a-1, a+1]$, the metric subspace $%
\mathbf{Y}$ of $\mathbb{R}$ is complete. The metric subspace $\mathbf{X}$ of $%
\mathbb{R}$, where $X=[a-1,a +1]$, is compact and connected. Clearly, the
set $Y$ is dense in $\mathbf{X}$. The identity function $f:\mathbf{Y}%
\rightarrow \mathbf{Y}$ is uniformly continuous. By Theorem \ref{t2.6}(i), $%
f$ does not extend continuously to a function from $Z=\{a\}\cup Y$ into $Y$.
Therefore, $f$ does not extend continuously over $X$. This, in view of
Proposition \ref{p1.1}, contradicts $\mathbf{UCE}_{cc}(\mathbb{R})$%
.\smallskip

(d) $\rightarrow $ (a) Assuming that $\mathbb{R}$ is sequential, we fix a metric
subspace $\mathbf{X}$ of $\mathbb{R}$, a complete metric subspace $\mathbf{Y}$ of 
$\mathbb{R}$, a dense subset $S$ of $\mathbf{X}$ and a uniformly continuous
function $f:\mathbf{S}\rightarrow \mathbf{Y}$. By our hypothesis and Theorem %
\ref{t2.3}(iii), $\mathbf{Y}$ is closed. So, by Corollary \ref{c3.2}, $f$
extends to a uniformly continuous function from $\mathbf{X}$ into $\mathbf{Y}
$.\smallskip

(ii) To prove (ii), we notice that it is known that $\omega -\mathbf{CAC}(\mathbb{R%
})$ holds true but $\mathbf{CAC}(\mathbb{R})$ fails in the
Fefferman-Levy Model $\mathcal{M}9$ of \cite{hr} (see, e.g., \cite{gg1} and Remark 4.59 of\cite{herl}). Therefore, in view of
Theorem \ref{t2.3}(iii)(b), $\mathbf{UCE}(\mathbb{R})$ also holds true in $%
\mathcal{M}9$.\medskip
\end{proof}

\begin{remark}
\label{r3.4} We can regard the interval $[0, 1]$ as the unique (up to an
equivalence) two-point Hausdorff compactification of $\mathbb{R}$. The unit
circle $S^1=\{x\in\mathbb{R}^2: \Vert x\Vert=1\}$, denoted by $%
\mathbf{S}^1$ as a topological subspace of $\mathbb{R}^2$, can be regarded
as the unique (up to an equivalence) one-point Hausdorff compactification of 
$\mathbb{R}$.
\end{remark}

The question which pops up at this moment is whether one can replace in $%
\mathbf{UCE}_{cc}(\mathbb{R})$ the two-point compactification $[0,1]$ with
the one-point compactification $\mathbf{S}^1$ of $\mathbb{R}$, i.e., if

\begin{itemize}
\item $\mathbf{UCE}_{1pc}(\mathbb{R})$ : For every complete metric subspace $%
\mathbf{Y}\ $of $\mathbb{R}^{2}$, if $D$ is a dense subset of the unit
circle $\mathbf{S}^1$ of $\mathbb{R}^{2}$ and $f:\mathbf{D}\rightarrow 
\mathbf{Y}$ is a uniformly continuous function, then there exists a unique
uniformly continuous extension $F: \mathbf{S}^1\to\mathbf{Y}$ of $f$.
\end{itemize}

\noindent is equivalent to $\mathbf{UCE}(\mathbb{R})$. We intend to show
that $\mathbf{UCE}_{1pc}(\mathbb{R})$ can be added to the list of
equivalents of Theorem \ref{t3.3} (i). We need the following new theorem
which leads to solutions of Problems 6.6 and 6.11 of \cite{OPWZ}:

\begin{theorem}
\label{t3.5} For each $m\in\mathbb{N}$, the space $\mathbb{R}^m$ is
sequential if and only if $\mathbb{R}$ is sequential.
\end{theorem}

\begin{proof}
Of course, if $\mathbb{R}^{m}$ is sequential, so is $\mathbb{R}$. Let us
assume that $\mathbb{R}$ is sequential. We shall prove by induction that $%
\mathbb{R}^{m}$ is sequential for each $m\in \mathbb{N}$. For $m\in \mathbb{N%
}$, let 
\begin{equation*}
\mathcal{B}(m)=\{B_{d_{e(m)}}(x,\frac{1}{n}):x\in \mathbb{Q}^{m},n\in 
\mathbb{N}\}.
\end{equation*}%
For $G\subseteq \mathbb{R}^{m}$ and each $i\in m$, let 
\begin{equation*}
G_{i}=\{t\in \mathbb{R}:\exists x\in G,x(i)=t\}.
\end{equation*}

We notice that, to prove that $\mathbb{R}^m$ is sequential, it suffices to
show that each sequentially closed subset of $\mathbb{R}^m$ is separable.

First, suppose that we have already proved that each bounded sequentially
closed subset of $\mathbb{R}^{m}$ is closed in $\mathbb{R}^{m}$, so compact.
Consider a sequentially closed subset $H$ of $\mathbb{R}^{m}$ and, for each $%
n\in \mathbb{N}$, put $K_{n}=\{x\in H:\Vert x\Vert\leq n\}$. The sets $K_{n}$ are
sequentially closed in $\mathbb{R}^{m}$ and bounded. By our hypothesis, the
sets $K_{n}$ are compact in $\mathbb{R}^{m}$. Using the base $\mathcal{B}(m)$%
, the compactness of every $K_{n}$ and a straightforward induction, we can
easily define a sequence $(\langle A_{n},\psi _{n}\rangle)_{n\in \mathbb{N}}$
such that each $A_{n}$ is a set dense in $K_{n}$, each $\psi _{n}$ is an
injection $\psi _{n}:A_{n}\rightarrow \mathbb{N}$. Then the set $%
A=\bigcup\limits_{n\in \mathbb{N}}A_{n}$ is countable and dense in $H$. If $a\in 
\overline{H}$, we can inductively define a sequence $(a_{n})_{n\in \mathbb{N}%
}$ of points of $A $ such that $\lim\limits_{n\to\infty} a_{n}=a$. Then $a\in H$, which shows
that $H$ is closed in $\mathbb{R}^{m}$.

We shall prove by induction that $\mathbb{R}^{m}$ is sequential. We have
already assumed that $\mathbb{R}^{1}=\mathbb{R}$ is sequential. Suppose that 
$m\in \mathbb{N}$ is such that $\mathbb{R}^{m}$ is sequential. Let $G$ be a
non-empty sequentially closed bounded subset of $\mathbb{R}^{m+1}$. Fix $%
i\in m+1$. We claim that $G_{i}$ is sequentially closed in $\mathbb{R}$. For 
$t\in G_{i}$, let $l(t)=\{x\in \mathbb{R}^{m+1}:x(i)=t\}$ and $L(t)=G\cap
l(t)$. Notice that the subspace $l(t)$ of $\mathbb{R}^{m+1}$ is closed, so
the set $L(t)$ is sequentially closed in $\mathbb{R}^{m+1}$. Since $l(t)$ is
homeomorphic with $\mathbb{R}^{m}$ and $\mathbb{R}^{m}$ is sequential, we
infer that $l(t)$ is sequential. This implies that $L(t)$ is closed in $l(t)$, so $L(t)$ is also closed in $%
\mathbb{R}^{m+1}$. Since $L(t)$ is also bounded, it is compact in $\mathbb{R}%
^{m+1}$. Now, fix a sequence $(x_{n})_{n\in \mathbb{N}}$ of points of $G_{i}$
such that $(x_{n})_{n\in \mathbb{N}}$ converges in $\mathbb{R}$ to $x$. For
each $n\in \mathbb{N}$, the set $L(x_{n})$ is a non-empty compact subset of $%
\mathbb{R}^{m+1}$. Therefore, using the base $\mathcal{B}(m+1)$, we can
inductively define a sequence $(y_{n})_{n\in \mathbb{N}}$ of points of $%
\mathbb{R}^{m+1}$ such that $y_{n}\in L(x_{n})$ for each $n\in \mathbb{N}$.
The sequence $(y_{n})_{n\in \mathbb{N}}$ is bounded, so it has a subsequence
which converges in $\mathbb{R}^{m+1}$ to a point $y$. For simplicity, we
assume that $(y_{n})_{n\in \mathbb{N}}$ converges to $y$. Then $y\in G$
because $G$ is sequentially closed. Since $y_{n}(i)=x_{n}$ for each $n\in 
\mathbb{N}$, we have $y(i)=x$. Hence $x\in G_{i}$. This proves that $G_{i}$
is sequentially closed in $\mathbb{R}$. Since $\mathbb{R}$ is sequential, $%
G_{i}$ is closed in $\mathbb{R}$. Hence, by Theorem 4.55 of \cite{herl},
there exists a countable set $D\subseteq G_{i}$ such that $D$ is dense in $%
G_{i}$. Let $D=\{t_{n}:n\in \mathbb{N}\}$. Using the base $\mathcal{B}(m+1)$
and the compactness of the non-empty sets $L(t_{n})$ for $n\in \mathbb{N}$,
we can inductively define a sequence $(\langle C_{n},\phi _{n}\rangle)_{n\in 
\mathbb{N}}$ such that each $C_{n}$ is a dense set in $L(t_{n})$ and each $%
\phi _{n}$ is an injective function from $C_{n}$ into $\mathbb{N}$. The set $%
C=\bigcup\limits_{n\in \mathbb{N}}C_{n}$ is countable and dense in $G$. Hence $G$
is closed in $\mathbb{R}^{m+1}$. By the mathematical induction, we have
proved that if $\mathbb{R}$ is sequential, so is $\mathbb{R}^{m}$ for each $%
m\in \mathbb{N}$.
\end{proof}

Now, we are in a position to prove the following theorem which is an
analogue for $\mathbb{R}^m$ of Theorem 4.55 of \cite{herl}:

\begin{theorem}
\label{t3.6} For each $m\in\mathbb{N}$, the following conditions are all
equivalent:

\begin{enumerate}
\item[(i)] $\mathbb{R}^m$ is sequential;

\item[(ii)] every complete metric subspace of $\mathbb{R}^m$ is closed in $%
\mathbb{R}^m$;

\item[(iii)] every complete and bounded metric subspace of $\mathbb{R}^m$ is
compact;

\item[(iv)] every sequentially compact subspace of $\mathbb{R}^m$ is compact;

\item[(v)] every complete metric subspace of $\mathbb{R}^m$ is separable;

\item[(vi)] every complete unbounded subspace of $\mathbb{R}^m$ contains an
unbounded sequence;

\item[(vii)] every non-empty countable family of non-empty complete metric
subspaces of $\mathbb{R}^m$ has a choice function;

\item[(viii)] every non-empty family of non-empty complete metric subspaces
of $\mathbb{R}^m$ has a choice function.
\end{enumerate}
\end{theorem}

\begin{proof}
First, we notice that a metric subspace $\mathbf{A}$ of $\mathbb{R}^m$ is complete
if and only if $A$ is sequentially closed in $\mathbb{R}^m$. That (i) and
(v) are equivalent has been shown in the proof to Theorem \ref{t3.5}. The
implications (i)$\rightarrow$ (ii)$\rightarrow$ (iii)$\rightarrow $(iv) and
(v)$\rightarrow$(vi), as well as (viii)$\rightarrow$(vii) are obvious.

Let $H\subseteq \mathbb{R}^{m}$ be sequentially closed in $\mathbb{R}^{m}$.
If $H$ is not closed in $\mathbb{R}^{m}$, then there exists $n_{0}\in 
\mathbb{N}$ such that the set $H_{n_{0}}=\{x\in H:\Vert x\Vert \leq n_{0}\}$
is not closed in $\mathbb{R}^{m}$. On the other hand, $H_{n_{0}}$ is
sequentially compact, so closed if (iv) holds. Hence (iv) implies (i).
Clearly, (vi) implies that every complete unbounded subspace of $\mathbb{R}$
contains an unbounded sequence. Hence, in view of Theorem 4.55 of \cite{herl}%
, (vi) implies that $\mathbb{R}$ is sequential. This, taken together with
Theorem \ref{t3.5}, proves that (vi) implies (i). In consequence, conditions
(i)-(vi) are all equivalent. To show that (i) implies (viii), let us
consider any non-empty family $\mathcal{F}$ of non-empty complete metric subspaces
of $\mathbb{R}^{m}$. For each $F\in \mathcal{F}$, let $n(F)$ be the least $%
n\in \mathbb{N}$ such that $\{x\in F:\Vert x\Vert \leq n\}\neq \emptyset $.
Assuming that (i) holds, we obtain that, for each $F\in \mathcal{F}$, the
set $K(F)=\{x\in F:\Vert x\Vert \leq n(F)\}$ is compact. It is known that
every non-empty family of non-empty compact subsets of $\mathbb{R}^{m}$ has
a choice function. Therefore, the family $\{K(F):F\in \mathcal{F}\}$ has a
choice function. Hence (i) implies (viii). To complete the proof, let us
deduce that (viii) implies (vi). Let $E$ be an unbounded complete metric subspace
of $\mathbb{R}^{m}$. Inductively, we can define a sequence $(n_{k})_{k\in
\omega }$ of natural numbers such that, for each $k\in \omega $, $%
n_{k}<n_{k+1}$ and $\emptyset \neq \{x\in E:n_{k}\leq \Vert x\Vert \leq
n_{k+1}\}$. For each $k\in \omega $, the metric subspace $E_{k}=\{x\in E:n_{k}\leq
\Vert x\Vert \leq n_{k+1}\}$ of $\mathbb{R}^{m}$ is complete. Assuming that
(viii) holds, we can choose $y\in \prod\limits_{k\in \omega }E_{k}$. Then $%
(y(k))_{k\in \omega }$ is an unbounded sequence of points of $E$. Hence
(viii) implies (vi).
\end{proof}

\begin{remark}
\label{r3.7} Let us notice that positive answers to the questions posed in
Problems 6.6 and 6.11 of \cite{OPWZ} follow immediately from Theorems \ref%
{t3.5} and \ref{t3.6}.
\end{remark}

\begin{remark}
\label{r3.d0} It follows directly from Theorem 4.54 of \cite{herl} and from
Exercise E.3 to Section 4.6 of \cite{herl} that, for each $n\in\mathbb{N}$, $%
\mathbf{CAC}(\mathbb{R})$ is equivalent with: $\mathbb{R}^n$ is
Fr\'echet-Urysohn.
\end{remark}

Let us leave a satisfactory answer to the following question for another
research project:

\begin{question}
If $\mathbb{R}$ is sequential, must $\mathbb{R}^{\omega}$ be sequential?
\end{question}

\begin{proposition}
\label{p3.8} Let $m\in\mathbb{N}$. If $\mathbb{R}$ is sequential, then every
closed subspace of $\mathbb{R}^m$ is sequential.
\end{proposition}

\begin{proof}
Assume that $\mathbb{R}$ is sequential. Let $\mathbf{X}$ be a closed
subspace of $\mathbb{R}^m$ and let $A\subseteq X$ be sequentially closed in $%
\mathbb{R}^m$. Since the set $X$ is closed in $\mathbb{R}^m$, the set $A$ is
sequentially closed in $\mathbb{R}^m$. In the light of Theorem \ref{t3.5},
the space $\mathbb{R}^m$ is sequential, so $A$ is closed in $\mathbb{R}^m$.
This implies that $A$ is closed in $\mathbf{X}$.
\end{proof}

\begin{proposition}
\label{p3.9} Let $m,n \in\mathbb{N}$ and let $\alpha\mathbb{R}^n$ be a
compactification of $\mathbb{R}^n$ such that $\alpha\mathbb{R}^n\subseteq%
\mathbb{R}^m$. Then $\mathbb{R}$ is sequential if and only if $\alpha\mathbb{%
R}^n$ is sequential.
\end{proposition}

\begin{proof}
If $\mathbb{R}$ is sequential, then $\alpha \mathbb{R}^{n}$ is sequential by
Proposition \ref{p3.8}. Now, assume that $\alpha \mathbb{R}^{n}$ is
sequential. Let $A$ be a sequentially closed subset of $\mathbb{R}^{n}$ and
let $C=\alpha \lbrack A]\cup (\alpha \mathbb{R}^{n}\setminus \alpha \lbrack 
\mathbb{R}^{n}])$. It is easy to notice that $C$ is sequentially closed in $%
\alpha \mathbb{R}^{n}$. Therefore, $C$ is closed in $\alpha \mathbb{R}^{n}$.
This implies that $A=\alpha ^{-1}[C]$ is closed in $\mathbb{R}^{n}$. Hence $%
\mathbb{R}^{n}$ is sequential and, in consequence, so is $\mathbb{R}$%
.\medskip
\end{proof}

Now, let us give a more general statement than $\mathbf{UCE}_{1pc}(\mathbb{R}%
)$. Namely, suppose that, for $m,n\in\mathbb{N}$, $\alpha\mathbb{R}^n$ is a
compactification of $\mathbb{R}^n$ such that $\alpha\mathbb{R}^n\subseteq%
\mathbb{R}^m$. Then we define $\mathbf{UCE}_{\alpha}(\mathbb{R}^n)$ as
follows:

\begin{itemize}
\item $\mathbf{UCE}_{\alpha}(\mathbb{R}^n)$ : For every complete metric
subspace $\mathbf{Y}$ of $\mathbb{R}^{m}$, if $D$ is a dense set in $\alpha%
\mathbb{R}^n $ and $f:\mathbf{D}\rightarrow \mathbf{Y}$ is a uniformly
continuous function, then there exists a uniformly continuous extension $%
F:\alpha\mathbb{R}^n\to\mathbf{Y}$ of $f$.
\end{itemize}

\begin{theorem}
\label{t.3.10} Assume that $m,n\in\mathbb{N}$ and that $\alpha\mathbb{R}^n$
is a compactification of $\mathbb{R}^n$ such that $\alpha\mathbb{R}%
^n\subseteq \mathbb{R}^m$. Then $\mathbf{UCE}_{\alpha}(\mathbb{R}^n)$ and $%
\mathbf{UCE}(\mathbb{R})$ are equivalent.
\end{theorem}

\begin{proof}
To show that $\mathbf{UCE}_{\alpha }(\mathbb{R}^n)$ implies $\mathbf{UCE}(%
\mathbb{R})$, assume that $\mathbf{UCE}_{\alpha }(\mathbb{R}^n)$ holds. In
view of Theorem \ref{t3.3}, it suffices to show that $\mathbb{R}$ is
sequential. To this end, we assume the contrary and fix, by Theorem \ref%
{t2.3}(iii)(d), a dense, sequentially closed, non-closed subset $A$ of $\mathbb{R}$.
Then the set $B=\{x\in\mathbb{R}^n: x(0)\in A\}$ is a sequentially closed,
non-closed dense subset of $\mathbb{R}^n$. The set $Y=\alpha \lbrack B]\cup
(\alpha \mathbb{R}^n\setminus \alpha \lbrack \mathbb{R}^n])$ is dense and
sequentially closed in $\alpha \mathbb{R}^n$. If $Y$ were closed in $\alpha 
\mathbb{R}^n$, the set $B$ would be closed in $\mathbb{R}^n$. Hence $Y$ is
not closed in $\alpha \mathbb{R}^n$. Let $a\in \overline{Y}\setminus Y$. The
metric subspace $\mathbf{Y}$ of $\mathbb{R}^{m}$ is complete, while, in view of
Theorem \ref{t2.6}(i), the identity function $f:\mathbf{Y}\rightarrow \mathbf{Y}
$ is not extendable to any continuous function $F:Y\cup \{a\}\rightarrow Y$.
This contradicts $\mathbf{UCE}_{\alpha }(\mathbb{R}^n)$. Hence $\mathbb{R}$
is sequential.

To show that $\mathbf{UCE}(\mathbb{R})$ implies $\mathbf{UCE}_{\alpha }(%
\mathbb{R}^n)$, we assume that $\mathbb{R}$ is sequential. Then, by Theorems %
\ref{t3.5} and \ref{t3.6}, every complete metric subspace of $\mathbb{R}^{m}$ is
closed, so $\mathbf{UCE}_{\alpha }(\mathbb{R}^n)$ follows from Corollary \ref%
{c3.2}.
\end{proof}

\begin{corollary}
\label{c3.11} $\mathbf{UCE}_{1pc}(\mathbb{R})$ and $\mathbf{UCE}(\mathbb{R})$
are equivalent.
\end{corollary}

Let us have a look at the following statement:

\begin{itemize}
\item $\mathbf{CAC}_D(\mathbb{R})$: Every family $\mathcal{A}=\{A_{n}: n\in 
\mathbb{N}\}$ of dense subsets of $\mathbb{R}$ has a partial choice function.
\end{itemize}

In contrast to the non-equivalence of $\mathbf{UCE}(\mathbb{R})$ and $%
\mathbf{CAC}(\mathbb{R})$, we are going to show that $\mathbf{DCC}(\mathbb{R}%
)$ is equivalent with $\mathbf{CAC}(\mathbb{R})$. The following theorem will
be applied:

\begin{theorem}
\label{td3.13} $\mathbf{CAC}(\mathbb{R})$ and $\mathbf{CAC}_{D}(\mathbb{R})$
are equivalent.
\end{theorem}

\begin{proof} Since it is obvious that $\mathbf{CAC}(\mathbb{R})$ implies $\mathbf{CAC}_D(\mathbb{R})$, we assume that $\mathbf{CAC}_D(\mathbb{R})$ holds and that $\mathcal{A}=\{A_{n}: n\in 
\mathbb{N}\}$ is a pairwise disjoint family of non-empty subsets of $\mathbb{R}$. We show that $\mathcal{A}
$ has a partial choice function. Let $((p_{n},q_{n}))_{n\in \mathbb{N}}$ be an
enumeration of the family of all open intervals with rational endpoints. Let 
$f: \mathbb{R}\rightarrow (0,1)$ be a homeomorphism. For each $n\in \mathbb{N}$, let  $g_n: \mathbb{R}\to\mathbb{R}$ be a linear function such that $g_n(0)=p_{n}$ and $g_{n}(1)=q_n$. For each $n\in\mathbb{N}$, we put $f_n=g_n\circ f$ and we consider the sets
\begin{equation*}
D_{n}=\bigcup \{f_{i}[A_{n}]: i\in \mathbb{N}\}\text{.}
\end{equation*}
It is easy to see that, for each $n\in\mathbb{N}$, $D_{n}$ is a dense
subset of $\mathbb{R}$.  By $\mathbf{CAC}_{D}(%
\mathbb{R})$, there exists an infinite subset $N$ of $\mathbb{N}$ such that $\mathcal{%
D}=\{D_{n}: n\in N\}$ has a choice function. Hence, there exists $h\in\prod\limits_{n\in N}D_n$. For each $n\in N$, let 
\begin{equation*}
k_{n}=\min \{i\in \mathbb{N}: h(n)\in f_{i}[A_{n}]\}\text{.}
\end{equation*}%
Since, for each $i\in\mathbb{N}$, the function $f_{i}$ is a bijection, it follows that $%
C=\{f_{k_{n}}^{-1}(h(n)): n\in N\}$ is a partial choice set of $\mathcal{A%
}$. 
\end{proof}

\begin{theorem}
\label{t3.12} For each $n\in\mathbb{N}$, the conditions $\mathbf{CAC}(%
\mathbb{R})$, $\mathbf{DCC}(\mathbb{R}^n)$ and $\mathbf{DCC}_c(\mathbb{R}^n)$
are all equivalent.
\end{theorem}

\begin{proof}
Let $n_0\in \mathbb{N}$ be given. We notice that $\mathbf{CAC}(\mathbb{R}^{n_0})$
and $\mathbf{CAC}(\mathbb{R})$ are equivalent because $\mathbb{R}$ and $%
\mathbb{R}^{n_0}$ are equipotent in $\mathbf{ZF}$. Therefore, that $\mathbf{CAC}(%
\mathbb{R})$ implies $\mathbf{DCC}(\mathbb{R}^{n})$ follows from Corollary %
\ref{c2.5}. It is straightforward that $\mathbf{DCC}_{c}(\mathbb{R})$
follows from $\mathbf{DCC}_{c}(\mathbb{R}^{n_0})$ and $\mathbf{DCC}(\mathbb{R}%
^{n_0})$ implies $\mathbf{DCC}_{c}(\mathbb{R}^{n_0})$. To complete the proof, it
suffices to show that if $\mathbf{CAC}(\mathbb{R})$ is false, then so is $%
\mathbf{DCC}_{c}(\mathbb{R})$.

Assume that $\mathbf{CAC}(\mathbb{R})$ is false. By Theorem \ref{td3.13}, we can fix a
family 
\begin{equation*}
\mathcal{D}=\{D_{p,q}:p,q\in \mathbb{Q}\cap (0,1],p<q\}
\end{equation*}%
of dense subsets of $\mathbb{R}$ having no partial choice function. For $p, q\in\mathbb{R}\cap (0, 1]$ with $p<q$,  we put $A_{p,q}=D_{p,q}\cap (p, q)$. The collection 
$$\mathcal{A}=\{ A_{p,q}: p,q\in \mathbb{Q}\cap (0,1],p<q\}$$
does not have a partial choice function, and the set $A=\bigcup\mathcal{A}$ is dense in the connected subspace $(0,1]$ of $\mathbb{R}$.

We claim that every Cauchy sequence of points of $A$ converges to some point
in $(0,1]$. To this end, it suffices to show that $0\notin \widetilde{A}$.
Assume the contrary and fix a strictly descending sequence $(x_{n})_{n\in
\omega}$ of points of $A$ converging to $0$. Via a straightforward
induction, we construct a subsequence $(x_{k_{n}})_{n\in \omega}$ of $%
(x_{n})_{n\in\omega}$ and a sequence of open intervals $((p_{n},q_{n}))_{n%
\in\omega}$ with rational endpoints in $(0,1]$ such that, for each $n\in 
\mathbb{N}$, $x_{k_{n+1}}<p_{n}$ and $x_{k_{n}}\in A_{p_{n},q_{n}}$.

First, we fix a well ordering on $\mathbb{Q}\times\mathbb{Q}$. For $n=0$, we
simply pick any $p_{0},q_{0}\in \mathbb{Q}\cap (0,1]$ such that $p_0<q_0$
and $x_{0}\in A_{p_{0},q_{0}}$. We put $k_{0}=0$.

Assume that $n\in\omega$ is such that the numbers $k_i$ and intervals $%
(p_{i},q_{i})$ have been defined for each $i\in n+1$ in such a way that, for
all $i,j\in n+1$ with $i<j$, $k_{i}<k_{j}$ and, moreover, the following
condition is satisfied: 
\begin{equation*}
x_{k_{i}}\in A_{p_{i},q_{i}}\text{ for each }i\in n+1,\text{
while }x_{k_{i+1}}<p_{i}\text{ for each }i\in n\text{.}
\end{equation*}%
Since $0<p_{n}$ and $(x_{n})_{n\in \mathbb{N}}$ converges to $0$, it follows
that infinitely many terms of $(x_{n})_{n\in \mathbb{N}}$ are strictly less
than $p_{n}$. Let $\langle p,q\rangle$ be the least member of $\mathbb{Q}%
\times \mathbb{Q}$ such that $0<p<q\leq 1$, $p<p_{n}$ and there exists $%
m\in\omega$ such that $k_n\in m$, while $x_{m}\in A_{p,q}$. Put 
\begin{equation*}
k_{n+1}=\min \{m\in\omega : k_{n}\in m\text{ and }x_{m}\in
A_{p,q}\},p_{n+1}=p,q_{n+1}=q
\end{equation*}%
to terminate the induction.

We notice that if $m,n\in\omega$ and $m\neq n$, then $\langle
p_{m},q_{m}\rangle\neq \langle p_{n},q_{n}\rangle$. Hence, if $h(\langle p_n, q_n\rangle)=x_{k_n}$ for each $n\in\omega$, then $h$ is a partial choice function of $\mathcal{A}$. This contradics our assumption on $\mathcal{D}$.  Thus, $(0,1]$
is densely complete as required.

Since if $\mathbf{CAC}(\mathbb{R})$ fails, then $(0,1]$ is a connected,
densely complete, not complete metric subspace of $\mathbb{R}$, it follows that $%
\mathbf{DCC}_c(\mathbb{R})$ implies $\mathbf{CAC}(\mathbb{R})$.
\end{proof}

\section{On $\mathbf{DCC}$ and $\mathbf{UCE}$}

\begin{theorem}
\label{t4.1} The following are equivalent:

\begin{enumerate}
\item[(i)] $\mathbf{DCC}$;

\item[(ii)] $\mathbf{CAC}$;

\item[(iii)] for every metric space $\mathbf{X}$ and every complete metric
subspace $\mathbf{Y}$ of $\mathbf{X}$, $\overline{Y}$ is complete;

\item[(iv)] for every metric space $\mathbf{X}$, every complete metric
subspace of $\mathbf{X}$ is closed in $\mathbf{X}$;

\item[(v)] functions between metric spaces are continuous if and only if
they are sequentially continuous;

\item[(vi)] (Cf. \cite{brun}.) real-valued functions on metric spaces are
continuous if and only if they are sequentially continuous;

\item[(vii)] every complete metric space is Fr\'{e}chet-Urysohn;

\item[(viii)] every complete metric space is sequential.
\end{enumerate}
\end{theorem}

\begin{proof}
That (ii) implies (i) follows from Corollary \ref{c2.5}.\smallskip\ 

(i)$\rightarrow $(ii) Assuming that $\mathbf{CAC}$ is false, we can consider
a pairwise disjoint family $\mathcal{A}=\{A_{n,m}:n,m\in \mathbb{N}\}$ of non-empty
sets having no partial choice set. Put $A=\bigcup \mathcal{A}$ and let $%
\{\infty _{n}:n\in \mathbb{N}\}$ be a set of pairwise distinct  elements, disjoint with $A$. For every $n\in 
\mathbb{N}$, let 
\begin{equation}
X_{n}=\bigcup \{A_{n,m}:m\in \mathbb{N}\}\cup \{\infty _{n}\},  \label{4}
\end{equation}%
and define a mapping $d_{n}:X_{n}\times X_{n}\rightarrow \mathbb{R}$ by
requiring:%
\begin{equation}
d_{n}(x,y)=d_{n}(y,x)=\left\{ 
\begin{array}{ll}
0\text{ if }x=y &  \\ 
\frac{1}{nm}\text{ if }x\in A_{n,m},y\in A_{n,s}\text{ and }m\leq s &  \\ 
\frac{1}{nm}\text{ if }x\in A_{n,m}\text{ and }y=\infty _{n} & 
\end{array}%
\right. \text{.}  \label{1}
\end{equation}%
It is a routine work to verify that for every $n\in \mathbb{N}$, $d_{n}$ is
a metric on $X_{n}$. Put $X=A\cup \{\infty _{n}:n\in \mathbb{N}\}$ and
define a function $d:X\times X\rightarrow \mathbb{R}$ by letting: 
\begin{equation*}
d(x,y)=d(y,x)=\left\{ 
\begin{array}{ll}
0\text{ if }x=y &  \\ 
\frac{1}{n}\text{ if }x\in X_{n},y\in X_{m}\text{ and }n<m &  \\ 
d_{n}(x,y)\text{ if }x,y\in X_{n} & 
\end{array}%
\right. \text{.}
\end{equation*}%
One can easily check that $d$ is a metric on $X$.\medskip 

\textbf{Claim ($\star $) }: $\overline{A}=X$ and every Cauchy sequence of $%
\mathbf{A}$ converges in $\mathbf{X}$.\medskip 

\textbf{Proof of Claim ($\star$)}. The first assertion of Claim ($\star$)
follows at once from the observation that, for every $n\in \mathbb{N}$, $%
\infty _{n}$ is an accumulation point of 
\begin{equation}
Y_{n}=\bigcup \{A_{n,m}:m\in \mathbb{N}\}  \label{5}
\end{equation}%
in the metric space $\mathbf{X}_{n}$. Hence, $\infty _{n}$ is an
accumulation point of $A$ in $\mathbf{X}$.

To see the second assertion, fix a Cauchy sequence $(x_{n})_{n\in \mathbb{N}%
} $ of $\mathbf{A}$. We show that $(x_{n})_{n\in \mathbb{N}}$ converges in $\mathbf{X}$. Assume the contrary and let $(x_{n})_{n\in \mathbb{N}}$
have no accumulation point in $\mathbf{X}$. Then $\{x_{n}:n\in \mathbb{N}\}$
is a countably infinite subset of $A$. By passing to a subsequence, if
necessary, we may assume that $(x_{n})_{n\in \mathbb{N}}$ is injective.
Since $\mathcal{A}$ has no partial choice set, it follows that there exist $%
k,m\in \mathbb{N}$ such that infinitely many terms of $(x_{n})_{n\in \mathbb{%
N}}$ belong to $A_{k,m}$. For our convenience let us assume that all the
terms of $(x_{n})_{n\in \mathbb{N}}$ belong to $A_{k,m}$. Then, $%
d(x_{t},x_{s})=\frac{1}{km}$ for all $t,s\in \mathbb{N},t\neq s$,
contradicting the fact that $(x_{n})_{n\in \mathbb{N}}$ is Cauchy, and
finishing the proof of Claim ($\star$).\medskip

It follows from Claim ($\star$) that $\mathbf{X}$ is densely complete.
However, $(\infty _{n})_{n\in \mathbb{N}}$ is easily seen to be a
non-convergent Cauchy sequence in $\mathbf{X}$. Hence $\mathbf{X}$ is not
complete, so $\mathbf{DCC}$ fails. This proves that (i) implies
(ii).\smallskip

The implications (ii)$\rightarrow$(iv)$\rightarrow$(iii), (ii)$\rightarrow$%
(v)$\rightarrow$(vi) and (ii)$\rightarrow$(vii)$\rightarrow$(viii) are
straightforward. To prove that (iii)$\rightarrow$(ii), (vi)$\rightarrow$(ii)
and (viii)$\rightarrow$(ii), let us assume that (ii) does not hold and let $%
\mathcal{A}$ be as in the proof of (i) $\rightarrow $ (ii). Let us use the
notation established in the proof of (i) $\rightarrow $ (ii).\smallskip

(iii)$\rightarrow$(ii) We claim that $\mathbf{A}$ is a complete and dense metric
subspace of $\mathbf{X}$ but $\mathbf{X}$, as in the proof of (i) $%
\rightarrow $ (ii), is not complete. This will contradict the (iii) and
terminate the proof of (iii) $\rightarrow $ (ii).

That $\overline{A}=X$ follows from Claim ($\star$). To see that $\mathbf{A}$
is complete, we fix a Cauchy sequence $(x_{n})_{n\in \mathbb{N}}$ of points
of $\mathbf{A}$. By Claim ($\star$), $(x_{n})_{n\in \mathbb{N}}$ converges
to some point $x\in \mathbf{X} $. If $x=\infty _{t}$ for some $t\in \mathbb{N%
},$ then almost all terms of $(x_{n})_{n\in \mathbb{N}}$ are included in $%
Y_{t}$, and for infinitely many $m\in \mathbb{N},$ $A_{t,m}\cap \{x_{v}:v\in 
\mathbb{N\}\neq \emptyset }$. Therefore, a partial choice function of $%
\mathcal{A}$ can be easily defined, contradicting our hypothesis for $%
\mathcal{A}$. Thus, $x\in A$ and $\mathbf{A}$ is complete.\smallskip
\smallskip

(vi) $\rightarrow $ (ii) We consider the sets $X_{1},Y_{1}$ defined by (\ref%
{4}) and (\ref{5}) respectively. Let $f:\mathbf{X}_{1}\rightarrow \mathbb{R}$
be the function given by: 
\begin{equation*}
f(x)=\left\{ 
\begin{array}{ll}
1\text{ if }x\in Y_{1} &  \\ 
0\text{ if }x=\infty _{1} & 
\end{array}%
\right. \text{.}
\end{equation*}%
Clearly, $f$ is a sequentially continuous (the only sequences of $\mathbf{X}%
_{1}$ converging to $\infty _{1}$ are those which are eventually equal to $%
\infty _{1}$), non-continuous (for $\varepsilon =1$ and every $\delta >0$, $%
|f(x)-f(\infty _{1})|\geq 1$ for every $x\in X_{1}$) function. This
contradicts (iii). Hence (iii) implies (ii).

(viii)$\rightarrow $(ii) Now, fix $n\in\mathbb{N}$ and consider the sets $%
X_{n},Y_{n}$ defined in the proof of (i)$\rightarrow $(ii). Working as in
the proof to Claim ($\star$), we can show that $\mathbf{X}_{n}$ is complete,
while $Y_n$ is sequentially closed in $\mathbf{X}_n$. Since $\infty_{n}\in 
\overline{Y_{n}}\setminus\widetilde{Y}_{n}$, we infer that $\mathbf{X}_n$ is not
sequential. Hence (viii) does not hold if (ii) is false.
\end{proof}

\begin{remark}
\label{r4.2} Since every Fr\'echet-Urysohn space is sequential in $\mathbf{ZF%
}$, while $\mathbf{CAC}$ implies that all metric spaces are
Fr\'echet-Urysohn, it follows from Theorem \ref{t4.1} that Theorem \ref{t2.3}%
(ii) holds. However, let us notice that our proof that conditions (ii),
(vii) and (viii) of Theorem \ref{t4.1} are equivalent is distinct than the
proof to Theorem \ref{t2.3}(ii) given in \cite{gg2}.
\end{remark}

\begin{theorem}
\label{t4.3} The following are equivalent:

\begin{enumerate}
\item[(i)] $\mathbf{CAC}$;

\item[(ii)] closed metric subspaces of densely complete metric spaces are
densely complete;

\item[(iii)] for every metric space $\mathbf{X}$, every sequential subspace $%
\mathbf{Y}$ of $\mathbf{X}$ and every $a\in\overline{Y}\setminus\widetilde{Y}
$, if $Z=Y\cup\{a\}$, then the subspace $\mathbf{Z}$ of $\mathbf{X}$ is
sequential.
\end{enumerate}
\end{theorem}

\begin{proof}
In view of Remark \ref{r4.2}, it is obvious that (i) implies (iii). 
By Theorem \ref{t4.1}, (i) implies densely
complete metric spaces are complete. Since closed metric subspaces of complete
metric spaces are complete, it is obvious that they are also densely
complete. Hence (i) implies (ii)

Now, let us assume that $\mathbf{CAC}$ does not hold. Fix a pairwise disjoint family $\mathcal{A}$
as in the proof of (i)$\rightarrow $(ii) of Theorem \ref{t4.1}. Let us also
use the notation established in the proof of Theorem \ref{t4.1}. Clearly,
the set $K=\{\infty _{n}:n\in \mathbb{N}\}$ is closed in $\mathbf{X}$, while
the metric subspace $\mathbf{K}$ of $\mathbf{X}$ fails to be densely complete. Moreover, the subspace $\mathbf{Y}_1$ of $\mathbf{X}_1$ is sequential, $\infty_1\in\overline{Y}_1\setminus\widetilde{Y_1}$, $X_1=Y_1\cup\{\infty_1\}$  but the space $\mathbf{X}_1$ is not sequential.
\end{proof}

The proof to the following Theorem \ref{t4.11} is another proof that $%
\mathbf{DCC}$ and $\mathbf{CAC}$ are equivalent. However, since the densely
complete, non-complete metric space $\mathbf{X}$ constructed in the proof to
Theorem \ref{t4.1} has been shown useful for other problems, we cannot omit
our first proof to the equivalence of (i) and (ii) of Theorem \ref{t4.1}

\begin{theorem}
\label{t4.11}  $\mathbf{DCC}_c$ and $\mathbf{CAC}$ are equivalent.
\end{theorem}

\begin{proof} Since $\mathbf{DCC}$ implies $\mathbf{DCC}_c$ and, in view of Corollary \ref{c2.5}(i), $\mathbf{CAC}$ implies $\mathbf{DCC}$, we infer that $\mathbf{DCC}_c$ follows from $\mathbf{CAC}$.

Now, suppose that $\mathbf{CAC}$ fails and fix a pairwise disjoint family $\mathcal{A}%
=\{A_{n}: n\in \mathbb{N}\}$ of non-empty sets having no partial choice
function. Fix a strictly increasing sequence $(a_{n})_{n\in\omega}$ of
points in $[0,1)$ converging to $1$. For each $n\in\omega$, let $r_{n}=\frac{a_{n+1}-a_{n}}{2},$ and $%
c_{n}=\frac{a_{n}+a_{n+1}}{2}$. In the complex plane $\mathbb{C}$, we consider the circles
$$C_n=\{ z\in\mathbb{C}: \vert z-c_n\vert=r_n\}$$
where $n\in\omega$. If $n\in\omega$ and $x\in A_n$, let $C_{n,x}=(C_n\setminus\{a_n, a_{n+1}\})\times\{x\}$.  Let 
\begin{equation*}
 Z=\{a_n: n\in\omega\}\cup\bigcup\{ C_{n,x}: n\in\omega \text{ and } x\in A_n\}. 
\end{equation*}
For every $n\in\omega$ and $a,b\in C_{n},$ let $\widehat{ab}$ denote the arc of the circle $C_{n}$ starting at $%
a$ and, moving counterclockwise, ending at $b$; moreover, let $l(ab)$ denote the
minimum of the lengths of the arcs $\widehat{ab},\widehat{ba}$, while $l_n=l(a_n a_{n+1})$.  Clearly, $%
l(ab)=l(ba) $. We define a metric $\rho $ on $Z$ by requiring: 
\begin{equation*}
\rho (\langle a,x\rangle,\langle b,x\rangle)=l(ab)\text{ and }\rho(a_i,\langle a,x\rangle)=\rho(\langle a,x\rangle, a_i)=l(a_i, a)
\end{equation*}%
in case $a,b\in C_{n}\setminus\{a_n, a_{n+1}\}$ and $x\in A_{n}$, while $n\in\omega$ and $i\in\{n, n+1\}$; 
\begin{equation}
\rho (\langle a,x\rangle, \langle b,y\rangle)=\min \{l(aa_{n+1})+l(ba_{n+1}),l(aa_{n})+l(ba_{n})\} \label{6}
\end{equation}%
in case $a,b\in C_n\setminus\{ a_n, a_{n+1}\}$ and $x,y\in A_{n}$, while $x\neq y$ and $n\in\omega$;
\begin{equation*}
\rho
(\langle a,x\rangle, \langle b,y\rangle)=\rho(\langle b,y\rangle, \langle a,x\rangle)=l(a_{n+1}a)+ l(bb_m)+\sum\limits_{i=n+1}^{m-1}l_i
\end{equation*}%
in case $a\in C_{n}\setminus\{a_n, a_{n+1}\}$, $b\in C_{m}\setminus\{a_m, a_{m+1}\}$ for some $n\in m\in\omega$, while $%
x\in A_{n},y\in A_{m}$;
\begin{equation*}
\rho(a_n, a_n)=0\text{ and }\rho(a_n, a_m)=\rho(a_m, a_n)=\sum\limits_{i=n}^{m-1} l_i
\end{equation*}
if $n\in m\in\omega$;
\begin{equation*}
\rho(a_n, \langle a, x\rangle)=\rho(\langle a, x\rangle, a_n)=l(aa_m)+\sum\limits_{i=n}^{m-1}l_i
\end{equation*}
when $n\in m\in\omega$, $a\in C_m\setminus\{a_m, a_{m+1}\}$ and $x\in A_m$;
\begin{equation*}
\rho(a_n, \langle a, x\rangle)=\rho(\langle a, x\rangle, a_n)=l(aa_{m+1})+\sum\limits_{i=m+1}^{n-1}l_i
\end{equation*}
when $m+1\in n\in\omega$, $a\in C_m\setminus\{a_m, a_{m+1}\}$ and $x\in A_m$.

Clearly, for each $n\in\omega$ and $x\in A_n$, the subspace $\{a_n, a_{n+1}\}\cup C_{n,x}$ of $\mathbf{Z}$ is pathwise connected. This implies that the metric space  $\mathbf{Z}=\langle Z, \rho\rangle$ is pathwise connected, so $\mathbf{Z}$ is connected. It is easy to see that the set 
\begin{equation*}
D=Z\setminus \{a_n: n\in\omega\}
\end{equation*}
is dense in $\mathbf{Z}$. Let us show that $\mathbf{Z}$ is densely complete with respect to $D$. To
this end, fix a Cauchy sequence $(d_{n})_{n\in\omega}$ of points of $D$.  Without loss of generality we may assume that $(d_{n})_{n\in\omega}$ is injective. Let
\begin{equation*}
B=\{d_n: n\in\omega\} \text{ and } \mathcal{C}=\{(C_m\setminus\{a_m, a_{m+1}\})\times A_m: m\in\omega\}.
\end{equation*}
If the set $\{C\in\mathcal{C}: C\cap B\neq\emptyset\}$ were infinite, using the countability of $B$, one could easily define a partial choice function of $\mathcal{A}$. Therefore, since $\mathcal{A}$ has no partial choice function, $B$ meets only finitely many members of $\mathcal{C}$. This implies that there exists $t\in\omega$ such that $B\subseteq\bigcup\{(C_m\setminus\{a_m, a_{m+1}\})\times A_m: m\in t+1\}$. There exists $k\in t+1$ such that the set $B\cap[(C_{k}\setminus\{a_k, a_{k+1}\})\times A_k]$
is infinite. For our convenience, let us assume
that $B\subseteq[(C_{k}\setminus\{a_k, a_{k+1}\})\times A_k]$. If there exists $x\in A_k$ such that  $B\cap C_{k,x}$ is infinite, then, by the completeness of the circle $\mathbf{C}_k$, the sequence $(d_{n})_{n\in \mathbb{N}}$ converges to some point of $C_{k,x}\cup\{a_k, a_{k+1}\}$.
Assume that, for each $x\in A_k$, the set $B\cap C_{k,x}$ is finite. We can replace, if necessary, the sequence $(d_n)_{n\in\omega}$ by  its subsequence $(d_{n_j})_ {j\in\omega}$ such that, for each $x\in A_k$, the set $\{d_{n_j}: j\in\omega\}\cap C_{k,x}$ consists of at most one point. Therefore, without loss of generality, we can assume that, for each $x\in A_k$, the set $B\cap C_{k,x}$ consists of at most one point. For each $n\in\omega$, let $d_{n}=\langle s_{n},x_{n}\rangle\in C_{k,x_{n}}$. We may assume that the sequence $(x_n)_{n\in\omega}$ is injective.  Since $%
(d_{n})_{n\in\omega}$ is Cauchy in $\mathbf{Z}$, it follows that there exists $n_{0}\in\omega$ such that if $n,m\in\omega$ and $n_0\subseteq n\cap m$, then
\begin{equation}
\rho (d_{n},d_{m})<\frac{\pi r_{k}}{3}\text{.}  \label{7}
\end{equation}%
Let $n\in\omega$ be such that $n_0\in n$. From \ref{6} and \ref{7}, it follows that
\begin{equation*}
\rho (d_{n},d_{n_{0}})=\min
\{l(s_{n}a_{k+1})+l(s_{n_{0}}a_{k+1}),l(s_{n}a_{k})+l(s_{n_{0}}a_{k})\}<\frac{\pi
r_{k}}{3}\text{.}
\end{equation*}%
Clearly, just one of $l(s_{n_{0}}a_{k+1})$ and $l(s_{n_{0}}a_{k})$ is less than $\frac{\pi
r_{k}}{3}$. Assume that $l(s_{n_{0}}a_{k})<\frac{\pi r_{k}}{3}$. Then $%
l(s_{n_{0}}a_{k+1})>\frac{2\pi r_{k}}{3}$ and $\rho
(d_{n},d_{n_{0}})=l(s_{n}a_{k})+l(s_{n_{0}}a_{k})<\frac{\pi r_{k}}{3}$. Hence $\lim\limits _{n\to\infty}l(s_{n}a_{k})=0$ and, in consequence, the sequence $(d_n)_{n\in\omega}$ converges in $\mathbf{Z}$ to $a_{k}$. This proves that $\mathbf{Z}$ is densely complete with respect  to $D$ as required. However, $(a_{n})_{n\in\omega}$ is easily seen to be a Cauchy
sequence in $\mathbf{Z}$ converging to no point of $Z$. This shows that $\mathbf{Z}$ is not complete and finishes the 
proof that $\mathbf{DCC}_c$ implies $\mathbf{CAC}$.
\end{proof}

\begin{definition}
Let $J$ be an at most countable non-empty set and let $\{ \mathbf{X}_j: j\in
J\}$ be a collection of metric spaces $\mathbf{X}_j=\langle X_j, d_j\rangle$
where $j\in J$. Let $X=\prod_{j\in J}X_j$.

\begin{enumerate}
\item[(i)] If $J$ is finite, let $d$ be the metric on $X$ defined as
follows: 
\begin{equation*}
d(x,y)=\max\{d_j(x(j), y(j)): j\in J\}
\end{equation*}
for all $x,y\in X$. Then $\mathbf{X}=\prod_{j\in J}\mathbf{X}_j$ denotes the
metric space $\langle X, d\rangle$ and $\mathbf{X}$ is called a \textit{%
finite product} of metric spaces.

\item[(ii)] If $J=\{j_{1},j_{2}\}$ where $j_{i}\neq j_{2}$, then $\mathbf{X}%
_{j_{1}}\times \mathbf{X}_{j_{2}}=\langle X_{j_{1}}\times X_{j_{2}},d\rangle 
$, where 
\begin{equation*}
d(\langle x_{1},y_{1}\rangle ,\langle x_{2},y_{2}\rangle )=\max
\{d_{j_{1}}(x_{1},x_{2}),d_{j_{2}}(y_{1},y_{2})\}
\end{equation*}%
for all $\langle x_{1},y_{1}\rangle ,\langle x_{2},y_{2}\rangle \in
X_{j_{1}}\times X_{j_{2}}$, is called the Cartesian product of the metric
spaces $\mathbf{X}_{j_{1}}$ and $\mathbf{X}_{j_{2}}$.

\item[(iii)] If $J=\{j_n: n\in\mathbb{N}\}$, we put $X_n= X_{j_n}$, $%
d_n=d_{j_n}$ and $\mathbf{X}_n=\mathbf{X}_{j_n}$ for each $n\in\mathbb{N}$.
Moreover, we put $X=\prod_{n\in\mathbb{N}} X_n$ and consider the metric $d$
on $X$ defined by (\ref{2}). Then $\mathbf{X}=\langle X, d\rangle$ is the
product $\prod_{n\in\mathbb{N}}\mathbf{X}_n$ of the collection $\{\mathbf{X}%
_j: j\in J\}$. For each $n\in\mathbb{N}$, let $\pi_n$ denote the projection $%
\pi_n: X\to X_n$ defined by: $\pi_n(x)=x(n)$ for each $x\in X$.
\end{enumerate}
\end{definition}

\begin{proposition}
\label{p4.5} Dense completeness is finitely productive, i.e., all finite
products of densely complete metric spaces are densely complete metric
spaces.
\end{proposition}

\begin{proof}
Let $J$ be a non-empty finite set and let $\{\mathbf{X}_{j}:j\in J\}$ be a
collection of metric spaces $\mathbf{X}_{j}=\langle X_{j},d_{j}\rangle $.
Suppose that, for each $j\in J$, the space $\mathbf{X}_{j}$ is densely
complete with respect to a set $D_{j}$. Then $\mathbf{X}=\prod_{j\in J}%
\mathbf{X}_{J}$ is densely complete with respect to the set $D=\prod_{j\in
J}D_{j}$.\medskip 
\end{proof}

In view of Proposition \ref{p4.5} and the countable productivity of
completeness for metric spaces in $\mathbf{ZF}$, see, e.g., \cite{kerem},
one may ask the following question:

\begin{question}
\label{q4.6} Is dense completeness in the class of metric spaces countably
productive in $\mathbf{ZF}$?
\end{question}

Our next two theorems give a partial answer to this question:

\begin{theorem}
\label{t4.7} Each of the following statements implies the one beneath it:

\begin{enumerate}
\item[(i)] $\mathbf{CMC}$.

\item[(ii)] For every family $\{\mathbf{X}_n:n\in \mathbb{N}\}$ of densely
complete metric spaces there exists a family $\{D_{n}:n\in \mathbb{N}\}$
such that, for every $n\in \mathbb{N}$, the space $\mathbf{X}_n$ is densely
complete with respect to $D_{n}$.

\item[(iii)] Every countable product of densely complete metric spaces is
densely complete.
\end{enumerate}
\end{theorem}

\begin{proof}
Let $\{\mathbf{X}_n : n\in \mathbb{N}\}$ be a collection of densely complete
metric spaces $\mathbf{X}_n=\langle X_n, d_n\rangle$ and let $\mathbf{X}%
=\prod_{n\in\mathbb{N}}\mathbf{X}_n$. If $\prod_{n\in\mathbb{N}}X_n=\emptyset
$, then $\mathbf{X}$ is trivially densely complete, so we may assume that $%
\prod_{n\in\mathbb{N}}X_n\neq\emptyset$.

(i)$\rightarrow$(ii) For each $n\in\mathbb{N}$, let $\mathcal{G}_n$ be the
collection of all sets $G$ such that $\mathbf{X}_n$ is densely complete with
respect to $G$. Assume that $\mathbf{CMC}$ holds. Then, by $\mathbf{CMC}$,
there exists a collection $\{\mathcal{D}_n: n\in\mathbb{N}\}$ such that, for
each $n\in\mathbb{N}$, $\mathcal{D}_n$ is a non-empty finite subcollection
of $\mathcal{G}_n$. For each $n\in\mathbb{N}$, we put $D_n=\bigcup\mathcal{D}%
_n$. It is straightforward that, for each $n\in\mathbb{N}$, the space $%
\mathbf{X}_n$ is densely complete with respect to $D_n$.\smallskip

(ii)$\rightarrow $(iii) Now, let us assume (ii). Fix a collection $\{D_n:
n\in\mathbb{N}\}$ such that, for each $n\in\mathbb{N}$, the space $\mathbf{X}%
_n$ is densely complete with respect to $D_n$. It is possible that $%
\prod_{n\in\mathbb{N}}D_n=\emptyset$. However, we have assumed that $%
\prod_{n\in\mathbb{N}}X_n\neq\emptyset$, so we can fix $a\in\prod_{n\in%
\mathbb{N}}X_n$. For $n\in\mathbb{N}$, let $H_n=D_n\cup\{a(n)\}$. Of course,
for each $n\in\mathbb{N}$, the space $\mathbf{X}_n$ is densely complete with
respect to $H_n$, while the set $H=\prod_{n\in\mathbb{N}}H_n$ is dense in $%
\mathbf{X}$. It can be easily deduced from the classical theory of metric
spaces that $\mathbf{X}$ is densely complete with respect to $H$.
\end{proof}

In \cite{kyr}, it has been shown that \textquotedblleft
sequentiality\textquotedblright\ of metric spaces is not countably
productive in $\mathbf{ZF}$. Next, we elaborate a little bit more on this
matter in connection with the notion of dense completeness.

\begin{theorem}
The following are equivalent: \label{t4.8}

\begin{enumerate}
\item[(i)] $\mathbf{CAC}$;

\item[(ii)] all countable products of densely complete metric spaces are
densely complete and sequential;

\item[(iii)] all countable products of discrete metric spaces are sequential;

\item[(iv)] all countable products of sequential metric spaces are
sequential.
\end{enumerate}
\end{theorem}

\begin{proof}
That (i) implies (ii) follows directly from Theorems \ref{t4.1} and \ref%
{t4.7}. It is known that (iv) follows from (i). It is obvious that (ii)
implies (iii) and (iv) implies (iii). To complete the proof, it suffices to
show that the negation of $\mathbf{CAC}$ implies the negation of (iii).

Assume that $\mathbf{CAC}$ does not hold and fix a pairwise disjoint family $\mathcal{A}%
=\{A_{n}:n\in \mathbb{N}\}$ of non-empty sets without a partial choice
function. Choose any set $\infty $ such that $\infty \neq \bigcup \mathcal{A}
$. For every $n\in \mathbb{N}$, let $X_{n}=A_{n}\cup \{\infty \}$ and let $%
d_{n}$ be the discrete metric on $X_{n}$. Put $\mathbf{X}_{n}=\langle
X_{n},d_{n}\rangle $ and $\mathbf{X}=\prod_{n\in \mathbb{N}}\mathbf{X}_{n}$.
Of course, all discrete spaces $\mathbf{X}_{n}$ are complete. We claim that $%
G=\bigcup_{n\in \mathbb{N}}\pi _{n}^{-1}[A_{n}]$ is sequentially closed in $%
\mathbf{X}$. To see this fix a sequence $(x_{n})_{n\in \mathbb{N}}$ of
points of $G$ which converges in $\mathbf{X}$ to a point $x$. Since $\mathcal{A}$ has no partial choice
function, it follows that there exists $n_{0}\in \mathbb{N}$ such that $%
(x_{n})_{n\in \mathbb{N}}\subseteq K$, where $K=\bigcup_{n\in n_{0}}\pi
_{n}^{-1}[A_{n}]$. Since $K$ is a closed subspace of $\mathbf{X}$, we see
that $x\in K\subseteq G$. Hence $G$ is sequentially closed in $\mathbf{X}$.
Let $a\in \ X$ be defined as follows: $a(n)=\infty$ for each $n\in 
\mathbb{N}$. Of course, $a\notin G$. Since, for every neighborhood $V$ of $a$
in $\mathbf{X}$, the set $V\cap G$ is non-empty, $a\in \overline{G}$. So, $G$
is sequentially closed but not closed in $\mathbf{X}$. This contradicts
(iii).\medskip 
\end{proof}

\begin{theorem}
\label{t4.9}  The following conditions are all equivalent:

\begin{enumerate}
\item[(i)] $\mathbf{CAC}(\mathbb{R})$;

\item[(ii)] for all non-empty metric subspaces $\mathbf{X}_{1},\mathbf{X}_{2}
$ of $\mathbb{R}$, if the product $\mathbf{X}_{1}\times \mathbf{X}_{2}$ is
densely complete, then both $\mathbf{X}_{1},\mathbf{X}_{2} $ are densely
complete;

\item[(iii)] for all non-empty topological subspaces $\mathbf{X}_{1},\mathbf{%
X}_{2}$ of $\mathbb{R}$, if the product $\mathbf{X}_{1}\times \mathbf{X}_{2}$
is densely completely metrizable, then both $\mathbf{X}_{1},\mathbf{X}_{2} $
are densely completely metrizable;

\item[(iv)] for each countably infinite, not complete metric subspace $%
\mathbf{C}$ of $\mathbb{R}$, the metric space $\mathbb{R}\times\mathbf{C}$
is not densely complete;

\item[(v)] the metric subspace $\mathbb{R}\times\{\frac{1}{n}: n\in\mathbb{N}%
\}$ of $\mathbb{R}^2$ is not densely complete;

\item[(vi)] the topological subspace $\mathbb{R}\times\mathbb{Q}$ of $%
\mathbb{R}^2$ is not densely completely metrizable.
\end{enumerate}
\end{theorem}

\begin{proof}
Since $\mathbf{CAC}(\mathbb{R})$ and $\mathbf{CAC}(\mathbb{R}^2)$ are equivalent, in the light of Theorem \ref{t2.4}(i), conditions (ii)-(vi) follow from (i). Now, assume that (i) is not satisfied. In view of Theorem \ref{td3.13}, we can fix a family $\mathcal{D}%
=\{D_{n}: n\in \mathbb{N}\}$ of dense subsets of $\mathbb{R}$ such that $\mathcal{D}$ does not have a partial choice function. Let $C=\{c_n: n\in\mathbb{N}\}$ be an infinitely countable subset of $\mathbb{R}$ and let $X=\mathbb{R}\times C$ be endowed with the Euclidean metric $\rho$. We may assume that $c_i\neq c_j$ for each pair $i,j$ of distinct numbers from $\mathbb{N}$.  Clearly, 
\begin{equation*}
D=\bigcup \{ D_{n}\times\{c_n\}: n\in \mathbb{N}\}
\end{equation*}
is a dense subset of $\mathbf{X}=\langle X, \rho\rangle$ such that $\{ D_{n}\times\{c_n\}: n\in 
\mathbb{N}\}$ has no partial choice function. To show that $\mathbf{X}$ is densely complete with respect to $D$, consider any Cauchy sequence $%
(z_{n})_{n\in \mathbb{N}}$ of points of $D$. Let $z_n=\langle x_n, y_n\rangle$ for each $n\in\mathbb{N}$. Then the set $\{y_n: n\in\mathbb{N}\}$ is finite because $\mathcal{D}$ does not have a partial choice function. Therefore, there exists $k\in\mathbb{N}$ such that $z_n\in \bigcup_{i=1}^k(D(i)\times\{c_i\})$ for each $n\in\mathbb{N}$. Hence $(z_n)_{n\in\mathbb{N}}$ is a Cauchy sequence of the complete metric subspace $\mathbb{R}\times\{c_1,\dots c_k\}$ of $\mathbb{R}^2$, so $(z_n)_{n\in\mathbb{N}}$ converges in $\mathbf{X}$. In consequence, if (i) does not hold, then all conditions (ii)-(vi) are false.
\end{proof}

\begin{remark}
\label{r4.10} The equivalence of (i) and (v) of Theorem \ref{t4.9} is
especially remarkable because $\mathbb{Q}$ is not densely completely
metrizable in $\mathbf{ZF}$, while $\mathbb{R}\times\mathbb{Q}$ is densely
completely metrizable in every model of $\mathbf{ZF}$ in which $\mathbf{CAC}(%
\mathbb{R})$ fails.
\end{remark}

\begin{theorem}
\label{t4.12} $\mathbf{UCE}$ and $\mathbf{CAC}$ are equivalent.
\end{theorem}

\begin{proof}
If $\mathbf{CAC}$ holds, then every complete metric space is Cantor complete; hence, by Lemma \ref{l3.1}, $\mathbf{UCE}$ follows from $\mathbf{CAC}$.  

In view of Theorem \ref{t4.1}, to prove that $\mathbf{CAC}$ follows from $\mathbf{UCE}$, it suffices to show that $\mathbf{UCE}$ implies
\textquotedblleft for every metric space $\mathbf{X}$, every complete
subspace $\mathbf{Y}$ of $\mathbf{X}$ is closed\textquotedblright . Assume
the contrary and fix a metric space $\mathbf{X}$ having a complete,
non-closed subset $Y$. Fix, $x\in \overline{Y}\backslash Y$. Clearly, the
identity function $f:\mathbf{Y\rightarrow Y}$ is uniformly continuous, and $%
\mathbf{Y}$ is dense in $\mathbf{M}$ where $M=\{x\}\cup Y$. By Theorem \ref{t2.6}(i),
$f$ cannot be extended to a continuous function on $\mathbf{M}$.
This contradicts $\mathbf{UCE}$ and finishes the proof.
\end{proof}

\begin{remark}
Clearly, $\mathbf{DCC}$ restricted to the class of compact metric spaces is
a theorem of $\mathbf{ZF}$. A compact metric space is always complete in $%
\mathbf{ZF}$. However, as Theorem \ref{t4.11} indicates, the restriction $%
\mathbf{DCC}_{c}$ of $\mathbf{DCC}$ to the class of all connected metric
spaces is not provable in $\mathbf{ZF}$. Similarly, Theorem \ref{t3.3} shows
that $\mathbf{UCE}_{cc}$ is unprovable in $\mathbf{ZF}$. Of course, $\mathbf{%
UCE}_c$ implies $\mathbf{UCE}_{cc}$.
\end{remark}

The following lemma will be applied in our proof that $\mathbf{UCE}_{cc}$
does not imply $\mathbf{CAC}$ in $\mathbf{ZF}^{0}$.

\begin{lemma}
\label{l4.14} Let $S$ be a subset of a metric space $\mathbf{X}=\langle X,
d\rangle$, let $\mathbf{Y}=\langle Y, \rho\rangle$ be a complete metric
space and let $f:\mathbf{S}\to\mathbf{Y}$ be a uniformly continuous mapping.
Then, for $Z=\widetilde{S}$, there exists a uniformly continuous mapping $F:%
\mathbf{Z}\to\mathbf{Y}$ such that $F(x)=f(x)$ for each $x\in S$.
\end{lemma}

\begin{proof}  Consider any  $x\in \widetilde{S}$. There exists a sequence $(x_n)_{n\in\omega}$ of points of $S$ which converges to $x$. Since $f$ is uniformly continuous, the sequence $(f(x_n))_{n\in\omega}$ is a Cauchy sequence in $\mathbf{Y}$, so it converges in $\mathbf{Y}$. Let $(x_n)_{n\in\omega}$ and $(x_n^{\star})_{n\in\omega}$ be sequences of points of $S$, both converging to $x$. Let $y$ be the limit in $\mathbf{Y}$ of $(f(x_n))_{n\in\omega}$ and let $y^{\star}$ be the limit in $\mathbf{Y}$ of $(f(x_n^{\star}))_{n\in\omega}$. Suppose that $y\neq y^{\star}$. Then $\varepsilon=\rho(y, y^{\star})>0$ and, by the uniform continuity of $f$, there exists $\delta>0$ such that $\rho(f(t), f(t^{\star}))<\frac{\varepsilon}{3}$ for all $t, t^{\star}\in S$ with $d(t, t^{\star})<\delta$.  There exists $k_1\in\omega$ such that $\rho(y, f(x_n))<\frac{\varepsilon}{3}$ and $\rho(y^{\star}, f(x_n^{\star}))<\frac{\varepsilon}{3}$ for each $n\in\omega\setminus k_1$.  There exists $k_2\in\omega$ such that if $n\in\omega\setminus k_2$,then  $d(x_n, x)<\frac{\delta}{2}$ and $d(x_n^{\star}, x)<\frac{\delta}{2}$. We notice that if $n\in\omega\setminus(k_1\cup k_2)$, then $d(x_n, x_n^{\star})\leq d(x_n, x)+d(x, x_n^{\star})<\delta$, so $\rho(f(x_n), f(x_n^{\star}))<\frac{\varepsilon}{3}$ and, therefore, $\rho(y, y^{\star})\leq\rho(y, f(x_n))+\rho(f(x_n), f(x_n^{\star}))+\rho(f(x_n)^{\star}, y^{\star})<\varepsilon$. This contradicts our assumption on $\varepsilon$. Hence $y=y^{\star}$. Let $y_x$ be the unique point of $Y$ such that, for every sequence $(x_n)_{n\in\omega}$ of points of $S$ converging to $x$, the sequence $(f(x_n))_{n\in\omega}$ is convergent in $\mathbf{Y}$ to $y_x$. We define a mapping $F:\mathbf{Z}\to \mathbf{Y}$ by putting $F(x)=y_x$ for each $x\in \widetilde{S}$. Of course, $F(x)=f(x)$ for each $x\in S$.  One can easily check in a standard way that $F$ is uniformly continuous. 
\end{proof}

\begin{theorem}
\label{t4.14}The conjunction of $\mathbf{BPI}$ and $\mathbf{CAC}(\mathbb{R})$
implies $\mathbf{UCE}_{c}$.\newline
In particular $\mathbf{UCE}_{c},$ and consequently $\mathbf{UCE}_{cc},$ does
not imply $\mathbf{CAC}$.
\end{theorem}

\begin{proof}
Let us assume $\mathbf{BPI}$ and $\mathbf{CAC}(\mathbb{R})$. Let $\mathbf{X}=\langle X, d\rangle$ be a compact metric space.  First, we show that $\mathbf{BPI}$ implies that 
the family $\mathcal{K}$ of all non-empty closed subsets of $\mathbf{X}$ has a choice function. Indeed, let $\infty$ be a set which does not belong to $X$ and,  for every $K\in\mathcal{K}$, let $\mathbf{X}_{K}$, where $X_K=K\cup \{\infty \}$, be the topological
sum of $\mathbf{K}$ and the trivial space $\{\infty \}$. Clearly, $K$ is a
closed subset of the compact space $\mathbf{X}_{K}$. Since 
\begin{equation*}
\mathcal{G}=\{\pi _{K}^{-1}(K): K\in \mathcal{K}\}
\end{equation*}
is a family of closed subsets of the Tychonoff product $\mathbf{Y}%
=\prod\limits_{K\in \mathcal{K}}\mathbf{X}_{K},$ and $\mathbf{Y}$ is compact
by $\mathbf{BPI}$ (see Theorem 4.70 in \cite{herl}), it follows that $\bigcap \mathcal{G}\neq \emptyset $.
Clearly, any element $f\in \bigcap \mathcal{G}$ is a choice function of $%
\mathcal{K}$. By Theorem \ref{t2.3}(iv)(b), $\mathbf{X}$ is separable. Hence, in view of Exercise E.3 to Section 4.6 of \cite{herl}, it follows from $\mathbf{CAC}(\mathbb{R})$ that every subspace of $\mathbf{X}$ is separable (see also \cite{kyr}). Now, suppose that $\mathbf{Y}=\langle Y, \rho\rangle$ is a complete metric space,  $S$ is a dense set in $\mathbf{X}$, and $f:\mathbf{S}\to\mathbf{Y}$ is uniformly continuous. Let $D$ be a countable dense set in $\mathbf{S}$. Since $X=\widetilde{D}$, it follows from Lemma \ref{l4.14} that there exists a uniformly continuous extension $F:\mathbf{X}\to\mathbf{Y}$ of $f$. Hence $\mathbf{UCE}_{c}$ holds in every model of $\mathbf{ZF}+[\mathbf{BPI}+\mathbf{CAC}(\mathbb{R})]$. 

The second assertion follows from the fact that in Mostowski's Linearly
Ordered Model $\mathcal{N}$3 in \cite{hr}, $\mathbf{BPI}$ and $\mathbf{CAC}(%
\mathbb{R})$ hold true, but $\mathbf{CAC}$ fails.
\end{proof}

\begin{remark}
As in the case of hedgehog metric spaces, it follows from the proof to
Theorem \ref{t4.11} that, in $\mathbf{ZFC}$, for every set $X$ such that $%
\vert X\vert\geq\vert\mathbb{R}\vert$, there exists a connected metric space 
$\mathbf{Z}$ with $|X|=|Z|$.
\end{remark}

We leave the following question unanswered:

\begin{question}
Does for every set $X$ equipotent with a subset of $\mathbb{R}$, there exist
in $\mathbf{ZF}$ a connected metric space $\mathbf{Z}$ equipotent with $X$?
\end{question}

\bigskip

\noindent \textsc{Kyriakos Keremedis}\newline
Department of Mathematics\newline
University of the Aegean\newline
Karlovassi, Samos 83200, Greece \newline
\emph{E-mail}: kker@aegean.gr\bigskip

\noindent\textsc{Eliza Wajch}

\noindent Institute of Mathematics and Physics

\noindent University of Natural Sciences and Humanities in Siedlce

\noindent ul. 3 Maja 54, 08-110 Siedlce, Poland

\noindent \emph{E-mail}: eliza.wajch@wp.pl

\end{document}